\documentclass[11pt,a4paper]{article}

\usepackage[utf8]{inputenc}
\usepackage[T1]{fontenc}
\usepackage{amsmath,amssymb,amsthm}
\usepackage{mathtools}
\usepackage{bm}
\usepackage[margin=1in]{geometry}
\usepackage{graphicx}
\usepackage{hyperref}
\usepackage{cleveref}
\usepackage{enumitem}
\usepackage{booktabs}
\usepackage{natbib}
\usepackage{xcolor}
\usepackage{algorithm} 
\usepackage{algorithmic}

\newtheorem{theorem}{Theorem}[section]

\newtheorem{proposition}[theorem]{Proposition}

\theoremstyle{definition}
\newtheorem{definition}[theorem]{Definition}

\newcommand{\R}{\mathbb{R}}
\newcommand{\C}{\mathbb{C}}

\newcommand{\BZ}{\mathcal{B}}
\newcommand{\kk}{\mathbf{k}}
\newcommand{\uu}{\mathbf{u}}
\newcommand{\xx}{\mathbf{x}}
\newcommand{\calM}{\mathcal{M}}
\newcommand{\calS}{\mathcal{S}}
\newcommand{\norm}[1]{\left\| #1 \right\|}

\newcommand{\dd}{\,\mathrm{d}}

\title{%
  On the Optimality of Reduced-Order Models for Band Structure Computations: \\[6pt]
  A Kolmogorov $n$-Width Perspective%
}

\author{
  Ankit Srivastava
}

\date{\today}

\begin{document}
\maketitle

\begin{abstract}
In this paper, we exploit the concept of Kolmogorov $n$-widths to establish
optimality benchmarks for reduced-order methods used in phononic, acoustic,
and photonic band structure calculations. The Bloch-transformed operators are
entire holomorphic functions of the wave vector~$\kk$, and by Kato's
analytic perturbation theory the eigenpairs inherit this holomorphy wherever
the spectral gap is positive. The Kolmogorov $n$-width of the solution
manifold therefore decays exponentially, at a rate controlled by the minimum
spectral gap between the band of interest and its neighbors. For clusters of
bands, we show that working with spectral projectors rather than individual
eigenvectors renders all internal crossings---avoided, symmetry-enforced, or
conical---irrelevant: only the gap separating the cluster from the remaining
spectrum matters. These results provide a sharp lower bound on the error of
any linear reduction method, against which existing approaches can be
measured. Numerical experiments on one- and two-dimensional problems confirm
the predicted exponential decay and demonstrate that a greedy algorithm
achieves near-optimal convergence. It also provides a principled justification for the choice of basis vectors in highly successful reduced-order models like RBME.
\end{abstract}

\medskip
\noindent\textbf{Keywords:} Kolmogorov $n$-width, band structure,
reduced-order models, Bloch eigenvalue problem, parametric holomorphy

\section{Introduction}

The propagation of waves through periodic media is governed by band structures.  In photonic crystals, periodic modulation of the dielectric permittivity gives rise to
photonic band gaps that can confine, guide, and manipulate electromagnetic
waves~\cite{Yablonovitch1987,John1987,Joannopoulos2008}. In phononic crystals, the
analogous modulation of elastic moduli and density produces band gaps for acoustic and elastic
waves, enabling vibration isolation, waveguiding, and exotic dynamic
phenomena~\cite{Sigalas1992,Kushwaha1993,MartinezSala1995,Hussein2014review,Srivastava2015review}.
In both cases, Bloch's theorem~\cite{Bloch1928,Brillouin1953} reduces the problem to a family of eigenvalue problems parametrized by the wave vector~$\mathbf{k}$
ranging over the Brillouin zone~$\mathcal{B}$.

Computing the band structure requires solving this eigenvalue problem at a large number of
wave vectors spanning~$\mathcal{B}$.  Each solve involves a system whose dimension~$N$ is determined by the spatial
discretization and can reach tens or hundreds of thousands of degrees of freedom. The total
computational cost---proportional to the product of the per-solve cost and the number of
$\mathbf{k}$-points---becomes a bottleneck, particularly when band structure evaluations are
embedded in outer loops such as topology optimization~\cite{LuYangGuestSrivastava2017}. This challenge has
motivated a sustained effort to develop efficient computational methods for band structure
evaluation in both the photonic and phononic settings.

\textit{The first generation} of methods focused on improving the efficiency and accuracy of the
\emph{full-order} eigenvalue solve at each $\mathbf{k}$-point.  The plane wave expansion (PWE)
method~\cite{Ho1990,Leung1990}, finite element methods~\cite{Dobson2000,ChinMokhtariSrivastavaSukumar2021}, and finite
difference time domain methods~\cite{Chan1995} each offer different trade-offs between
generality, accuracy, and computational cost.  Variational formulations using both
displacement and stress fields---the mixed variational
approach~\cite{NematNasser1972,NematNasser1975,SrivastavaNematNasser2014,LuSrivastava2016,LuSrivastava2017}---achieve
faster convergence than displacement-only formulations for unit cells with discontinuous
material properties. Dedicated eigenvalue solvers 
and iterative methods~\cite{Johnson2001} further reduce the per-solve cost. While these methods make individual solves fast, the
fundamental cost structure remains: the eigenvalue problem must be solved independently at
every $\mathbf{k}$-point, with no information reused across the Brillouin zone.

\textit{A second generation} of methods addresses this limitation through \emph{reduced-order
modeling}: the idea that eigenvectors computed at a few selected wave vectors can be assembled
into a compact basis that accurately represents the solution across the entire Brillouin zone.
All such methods share the structure described in Section~2.4 of this paper: a small set of
basis vectors is chosen, the eigenvalue problem is projected onto the resulting low-dimensional
subspace, and the reduced problem is solved inexpensively at each~$\mathbf{k}$-point. The Reduced Bloch Mode Expansion (RBME) method of
Hussein~\cite{Hussein2009} constructs the reduced basis from Bloch eigenvectors computed at the
high-symmetry points of the irreducible Brillouin zone. This Bloch-mode selection strategy builds on
earlier multiscale methods~\cite{HusseinHulbert2006} and on analogous ideas in electronic
structure calculations~\cite{Shirley1986}. Being a secondary expansion applied on top of any
primary discretization, RBME preserves the favorable properties of the underlying method
while reducing computation time by up to two orders of magnitude.

An alternative family of reduction methods adapts the Craig--Bampton component mode
synthesis (CMS) framework~\cite{CraigBampton1968} to the Bloch eigenvalue problem.  In this
approach, the unit cell is partitioned into interior and boundary degrees of freedom; the
interior is reduced via fixed-interface normal modes, and constraint modes account for the
boundary.  Krattiger and Hussein introduced this as Bloch Mode Synthesis
(BMS)~\cite{KrattigerHussein2014} and subsequently generalized it with residual-mode
enhancement and local boundary reduction, achieving speedups of one to three orders of
magnitude~\cite{KrattigerHussein2018}. The methodology has since been extended in several
directions: Palermo and Marzani developed Extended BMS for the computation of complex band
structures (the $\mathbf{k}(\omega)$ problem)~\cite{PalermoMarzani2016}; Xi and Zheng
incorporated algebraic condensation to improve
efficiency~\cite{XiZheng2021}; and Aladwani, Nouh, and Hussein generalized BMS to
state-space form for non-classically damped phononic
materials~\cite{AladwaniNouhHussein2022}. In photonics, multipoint model-order reduction
schemes based on the finite element method have been developed for efficient computation of
photonic crystal band diagrams~\cite{ScheiberSchultschikBiroDyczijEdlinger2011}.

These methods demonstrate impressive computational speedups with negligible loss of accuracy
across a range of problems.  Yet a fundamental question has remained unaddressed: \emph{how
close are these methods to the best possible approximation achievable by any $n$-dimensional
linear subspace?} This question has a precise mathematical formulation through the concept of the Kolmogorov
$n$-width~\cite{Kolmogorov1936,Pinkus1985}, which measures the best-possible worst-case
approximation error achievable by \emph{any} $n$-dimensional subspace.  In the context of
parametric partial differential equations, a body of work in approximation theory has
established that when the parameter-to-solution map is holomorphic---extendable as an analytic
function into a complex neighborhood of the parameter domain---the $n$-width decays
exponentially with~$n$, and the decay rate is controlled by the size of the holomorphic
extension~\cite{CohenDeVoreSchwab2011,CohenDeVore2015,BachmayrCohen2017}. These results
provide theoretical justification for the success of reduced basis and POD methods for
parametric elliptic problems~\cite{RozzaHuynhPatera2008,DeVore2017}. The greedy algorithm
for reduced basis construction has been shown to be rate-optimal: its error decays at the same
rate as the $n$-width~\cite{BinevCohenDahmenDeVorePetrovaWojtaszczyk2011}---providing a practical recipe for constructing a near optimal basis. Conversely, for
problems where the solution depends non-smoothly on the parameter---such as transport problems
with moving discontinuities---the $n$-width decays only algebraically, fundamentally limiting
the efficiency of any linear reduction method~\cite{GreifUrban2019,OhlbergerRave2016}.

While the connection between parametric regularity and reducibility is well established for
source problems governed by elliptic PDEs, it has not been systematically applied to the Bloch
eigenvalue problem.  The present paper bridges this gap.  We show that the Bloch eigenvalue
problem is a particularly clean instance of the parametric holomorphy framework, and the
eigenpairs inherit holomorphic dependence on~$\mathbf{k}$ wherever the spectral gap is
positive, by Kato's analytic perturbation theory~\cite{Kato1995}. The Kolmogorov $n$-width of
the solution manifold therefore decays exponentially, with a rate controlled by the minimum
spectral gap between the band of interest and its neighbors. 

The paper is organized as follows. Section~2 formulates the Bloch eigenvalue problem,
establishes the affine structure in~$\mathbf{k}$, and defines the solution manifold and the
reduced-order approximation framework. Section~3 introduces the Kolmogorov $n$-width,
discusses its role as an optimality benchmark, and reviews the key results connecting
parametric holomorphy to exponential $n$-width decay.  Section~4 applies this framework to the
Bloch problem: we establish holomorphy of the eigenpairs, derive $n$-width bounds for isolated
bands and for multi-band manifolds via spectral projectors, and discuss the interplay between
domain, band selection, and resolution type.  Section~5 presents numerical explorations for a
one-dimensional phononic crystal, comparing the SVD-optimal subspace, an oracle greedy
algorithm, and a practical residual-based greedy against the theoretical predictions. Section 6 presents numerical results for a 2-dimensional periodic composite.

\subsection{Problem formulation}
\label{sec:bloch:formulation}

Consider a periodic medium occupying all of $\R^d$ ($d = 1, 2, 3$) with a
unit cell $\Omega \subset \R^d$ and lattice vectors
$\{\mathbf{a}_1, \ldots, \mathbf{a}_d\}$. The medium may support elastic,
acoustic, or electromagnetic waves, depending on the physical context. In each case, the material properties are periodic with respect to the
lattice, and Bloch's theorem~\cite{Bloch1928,Brillouin1953} reduces the wave
equation on the infinite domain to a family of eigenvalue problems on the
unit cell, parametrized by the wave vector $\kk$ ranging over the first
Brillouin zone $\BZ$.

For elastic waves, the Bloch-transformed equation takes the form of a
generalized eigenvalue problem for the periodic part of the displacement
field. For scalar acoustic waves, the same structure arises with the
stiffness operator built from $\kappa^{-1}(\xx)$ and the mass operator from
$\rho(\xx)$. For electromagnetic waves, the standard formulation casts the
problem as an eigenvalue equation for the magnetic field
$\mathbf{H}$~\cite{Joannopoulos2008}, with the curl-curl operator weighted by
$\varepsilon^{-1}(\xx)$ playing the role of the stiffness operator.

After spatial discretization (by finite elements, plane waves, or other
methods), all three settings yield a $\kk$-dependent generalized matrix
eigenvalue problem on the unit cell:
\begin{equation}
  \label{eq:bloch-evp}
  K(\kk)\,\uu = \omega^2\,M(\kk)\,\uu, \qquad \kk \in \BZ,
\end{equation}
where $K(\kk)$ and $M(\kk)$ are Hermitian positive-semidefinite matrices
whose $\kk$-dependence arises from the Bloch-periodic boundary conditions
imposed on the unit cell (as detailed in \Cref{sec:bloch:affine}), and
$\omega^2$ and $\uu$ are the eigenvalue and eigenvector, respectively. The
dimension $N$ of the system is determined by the spatial discretization and
the number of field components: $N$ can range from hundreds for simple
two-dimensional scalar problems to hundreds of thousands for
three-dimensional vector problems with complex unit cell geometries.

The entire analysis that follows---the affine decomposition, the holomorphy
of the eigenpairs, and the $n$-width bounds---depends only on the abstract
structure~\eqref{eq:bloch-evp} and applies uniformly to phononic, acoustic,
and photonic band structure calculations. For concreteness, we use the
language of phononic crystals (stiffness, mass, displacement) throughout, but
all results hold verbatim for the other settings.

\subsection{Affine structure in $\kk$}
\label{sec:bloch:affine}

A key structural property of the Bloch problem is that the operators $K(\kk)$
and $M(\kk)$ depend on $\kk$ through a finite number of trigonometric
functions. This structure arises naturally from the imposition of Bloch-periodic boundary
conditions on the unit cell. It can be shown (see \Cref{app:1d-example} in the Appendix for an example) that the stiffness matrix is a Laurent polynomial in the phase factors:
\begin{equation}
  \label{eq:affine-decomp}
  K(\kk) = K_0 + \sum_{m=1}^{Q} f_m(\kk)\,K_m,
\end{equation}
where the $K_m$ are $\kk$-independent matrices determined by the element
stiffness contributions, and each coefficient function $f_m(\kk)$ is a
monomial of the form
\begin{equation}
  \label{eq:phase-monomials}
  f_m(\kk) = \prod_{j=1}^{d} e^{i\,\alpha_{mj}\,\kk\cdot\mathbf{a}_j},
  \qquad \alpha_{mj} \in \{-1, 0, 1\},
\end{equation}
or a sum of such monomials. The number of terms $Q$ is finite and determined
by the connectivity of the finite element mesh across the unit cell
boundaries; it does not grow with mesh refinement. The mass matrix $M(\kk)$
admits the same type of decomposition. Since the coefficient functions $f_m(\kk)$ are built from exponentials
$e^{i\kk\cdot\mathbf{a}_j}$, each $f_m$ extends to an entire function of
$\kk \in \C^d$---it is holomorphic on all of $\C^d$ with no singularities.
Consequently, the operator-valued maps $\kk \mapsto K(\kk)$ and
$\kk \mapsto M(\kk)$ are entire holomorphic families.

\subsection{Solution manifold and reduced order strategy}
\label{sec:bloch:manifold}

For each band index $j$, the $j$-th eigenvector defines a mapping
$\kk \mapsto \uu_j(\kk)$ from the Brillouin zone into the solution space $V$.
The \emph{solution manifold} for the $j$-th band is
\begin{equation}
  \label{eq:sol-manifold}
  \calM_j = \bigl\{ \uu_j(\kk) \in V : \kk \in \BZ \bigr\}.
\end{equation}

The central computational task in phononic band structure analysis is the
solution of the eigenvalue problem \eqref{eq:bloch-evp} at a large number of
wave vectors $\kk$ across the Brillouin zone $\BZ$. If the discretized
problem has dimension $N$---as determined by the finite element mesh or other
spatial discretization---then each eigenvalue solve carries a cost that grows
rapidly with $N$, and this cost is incurred at every $\kk$-point. 

All efficient methods for alleviating this cost share, at a fundamental level,
the same strategy. One selects, by whatever means, a set of $n$ vectors
$\{\boldsymbol{\phi}_1, \ldots, \boldsymbol{\phi}_n\}$ in the full
$N$-dimensional space $V$, with $n \ll N$, and forms the subspace
$V_n = \mathrm{span}\{\boldsymbol{\phi}_1, \ldots, \boldsymbol{\phi}_n\}$.
The original eigenvalue problem \eqref{eq:bloch-evp} is then replaced by its
projection onto $V_n$: assembling the reduced matrices
\begin{equation}
  \label{eq:reduced-evp}
  \widetilde{K}(\kk) = \Phi^\top K(\kk)\,\Phi, \qquad
  \widetilde{M}(\kk) = \Phi^\top M(\kk)\,\Phi,
\end{equation}
where $\Phi = [\boldsymbol{\phi}_1 \;\cdots\; \boldsymbol{\phi}_n] \in
\R^{N \times n}$ collects the basis vectors, and solving the reduced problem
\begin{equation}
  \label{eq:reduced-evp-solve}
  \widetilde{K}(\kk)\,\widetilde{\uu} = \widetilde{\omega}^2\,
  \widetilde{M}(\kk)\,\widetilde{\uu}, \qquad
  \widetilde{\uu} \in \R^n.
\end{equation}
Since $n \ll N$, the reduced eigenvalue problem \eqref{eq:reduced-evp-solve}
is inexpensive to solve. The accuracy of this approximation depends entirely on how well the subspace
$V_n$ captures the true eigenvectors $\uu_j(\kk)$ across the Brillouin zone or how closely
$V_n$ approximates the solution manifold $\calM_j$---or, if multiple bands
are of interest, the union $\calM = \bigcup_j \calM_j$. If $V_n$ is chosen
well, the projected eigenvectors $\Phi\,\widetilde{\uu}$ will be close to
the true eigenvectors $\uu_j(\kk)$ for all $\kk \in \BZ$, and the
approximate eigenvalues will be close to the true eigenvalues.

This shared structure raises a natural and fundamental
question: \emph{given a fixed dimension $n$, what is the best possible choice
of $V_n$?} Answering
this question requires a framework for measuring the intrinsic
approximability of the solution manifold $\calM$ by linear subspaces of
prescribed dimension. Such a framework is provided by the theory of
Kolmogorov $n$-widths.

\section{Kolmogorov $n$-Widths}
\label{sec:nwidth}

Consider a parametric family of problems---indexed by a parameter $\mu$ ranging over
some set $\mathcal{P}$---and for each value of $\mu$, there is a
corresponding solution $\uu(\mu)$ living in a high-dimensional space $V$ (for
instance, a finite element space of dimension $N$). The collection of all such
solutions,
\begin{equation}
  \label{eq:sol-manifold-generic}
  \calM = \bigl\{\uu(\mu) \in V : \mu \in \mathcal{P}\bigr\},
\end{equation}
is the \emph{solution manifold}. Although $V$ may be high dimensional, $\calM$ is typically a
low-dimensional object: it is parametrized by $\mu$, which often lives in a
space of dimension $p \ll N$. The practical question is whether one can find a subspace $V_n \subset V$ of
small dimension $n \ll N$ such that every element of $\calM$ is well
approximated by its projection onto $V_n$. If such a subspace exists, one can
replace the original $N$-dimensional computation with an $n$-dimensional one
at greatly reduced cost. This is the underlying idea behind all reduced-order
modeling and modal reduction methods.

\begin{definition}[Kolmogorov $n$-width]
  \label{def:nwidth}
  Let $\calM$ be a compact subset of a normed space $(V, \norm{\cdot})$. The
  \emph{Kolmogorov $n$-width} of $\calM$ in $V$ is
  \begin{equation}
    \label{eq:nwidth-def}
    d_n(\calM, V) = \inf_{\substack{V_n \subset V \\ \dim(V_n) = n}}
    \;\sup_{\uu \in \calM}
    \;\inf_{\mathbf{v} \in V_n} \norm{\uu - \mathbf{v}}_V.
  \end{equation}
\end{definition}

\noindent The inner infimum, $\inf_{\mathbf{v} \in V_n}\norm{\uu - \mathbf{v}}$,
    measures how well a single solution $\uu$ is approximated by the subspace
    $V_n$. In a Hilbert space, this is simply the distance from $\uu$ to $V_n$,
    realized by the orthogonal projection. The supremum, $\sup_{\uu \in \calM}$, selects the worst-case solution
    over the entire manifold. The outer infimum, $\inf_{V_n}$, optimizes over all possible
    $n$-dimensional subspaces.

The quantity $d_n(\calM, V)$ is therefore an intrinsic property of the solution
manifold $\calM$ itself. It depends on the physics of the problem but not on any particular computational method used to construct a
reduced basis. The sequence $\{d_n\}_{n \geq 0}$ is non-increasing: adding a dimension to
the approximation subspace can only help.

The rate at which $d_n$
decreases with $n$ is the central object of study.

\begin{figure}[htp]
    \centering
    \includegraphics[width=\textwidth]{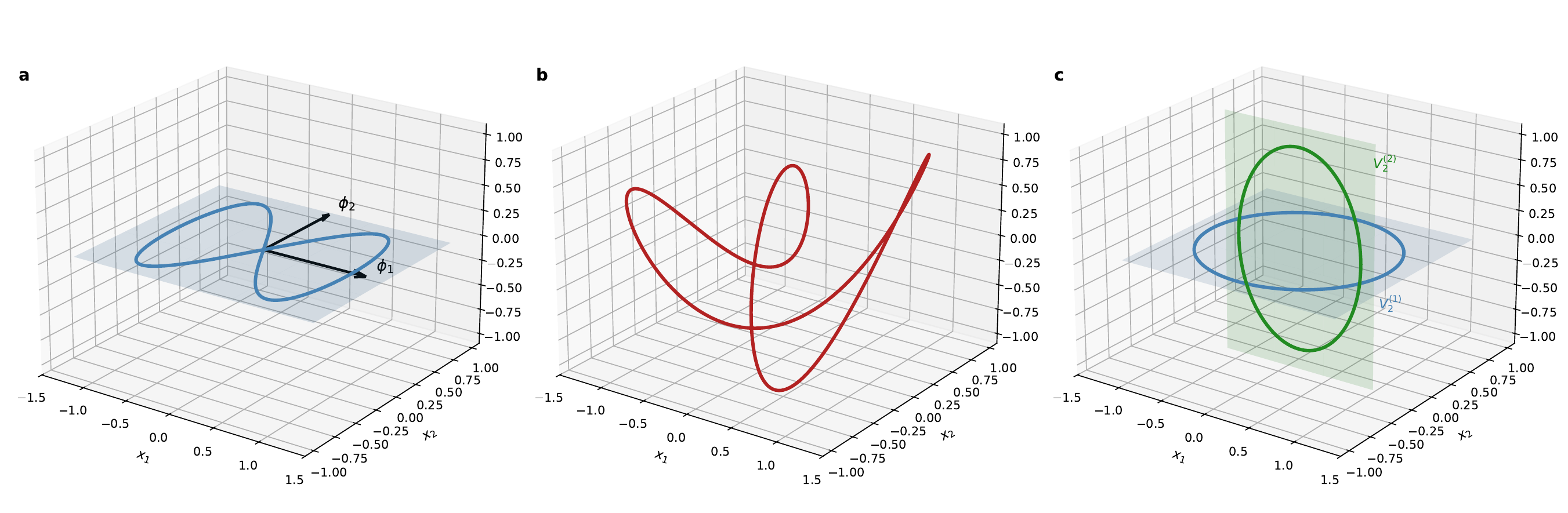}
    \caption{One-dimensional parametric curves (solution manifolds) embedded
    in a three-dimensional ambient space $V = \mathbb{R}^3$.
    \textbf{(a)}~A curve lying entirely within a two-dimensional subspace
    $V_2$: two basis vectors $\phi_1, \phi_2$ suffice to represent every point
    on $\mathcal{M}$ exactly.
    \textbf{(b)}~A curve winding through all three dimensions.
    \textbf{(c)}~Two separate one-dimensional curves, each confined to its own
    two-dimensional subspace which together span all three dimensions.}
    \label{fig:manifold}
\end{figure}

As an explanatory example, Fig.~\ref{fig:manifold} illustrates an example manifold geometry
in a three-dimensional ambient space $V = \mathbb{R}^3$. In panel~(a), the
solution manifold is a one-dimensional curve that happens to lie within a
two-dimensional subspace $V_2$; two basis vectors $\phi_1, \phi_2$ suffice to
represent every point on $\mathcal{M}$ exactly, so $d_2(\mathcal{M}, V) = 0$.
In panel~(b), the manifold is again a one-dimensional curve but it winds through all three dimensions, so no
two-dimensional subspace captures it well. The key point is that the intrinsic dimension
of $\mathcal{M}$ (here, one) does not determine the dimension $n$ of the
subspace needed: it is the \emph{geometry} of how $\mathcal{M}$ curves
through $V$ that controls how efficiently $\mathcal{M}$ can be
approximated by a linear subspace.

Panel~(c) illustrates the multi-curve analog. Two separate one-dimensional
curves are shown, each confined to its own two-dimensional subspace; their union is still
intrinsically one-dimensional, yet no single two-dimensional subspace captures
both curves simultaneously---together they span all three dimensions. This is
directly analogous to the multi-band setting of the Bloch eigenvalue problem.
For an $m$-dimensional crystal, the solution manifold of a single band is
$m$-dimensional, parametrized by the wave vector $\kk \in \mathcal{BZ}
\subset \mathbb{R}^m$. When $J$ bands are targeted simultaneously, the union
$\mathcal{M} = \bigcup_{j=1}^{J} \mathcal{M}_j$ remains $m$-dimensional as a
geometric object---it consists of $J$ sheets each parametrized by the same
$\kk$---but a single subspace $V_n$ must accommodate all $J$ sheets at once.

\subsection{The $n$-width as an optimality benchmark}
\label{sec:nwidth:benchmark}

\begin{figure}[htbp]
    \centering
    \includegraphics[width=\textwidth]{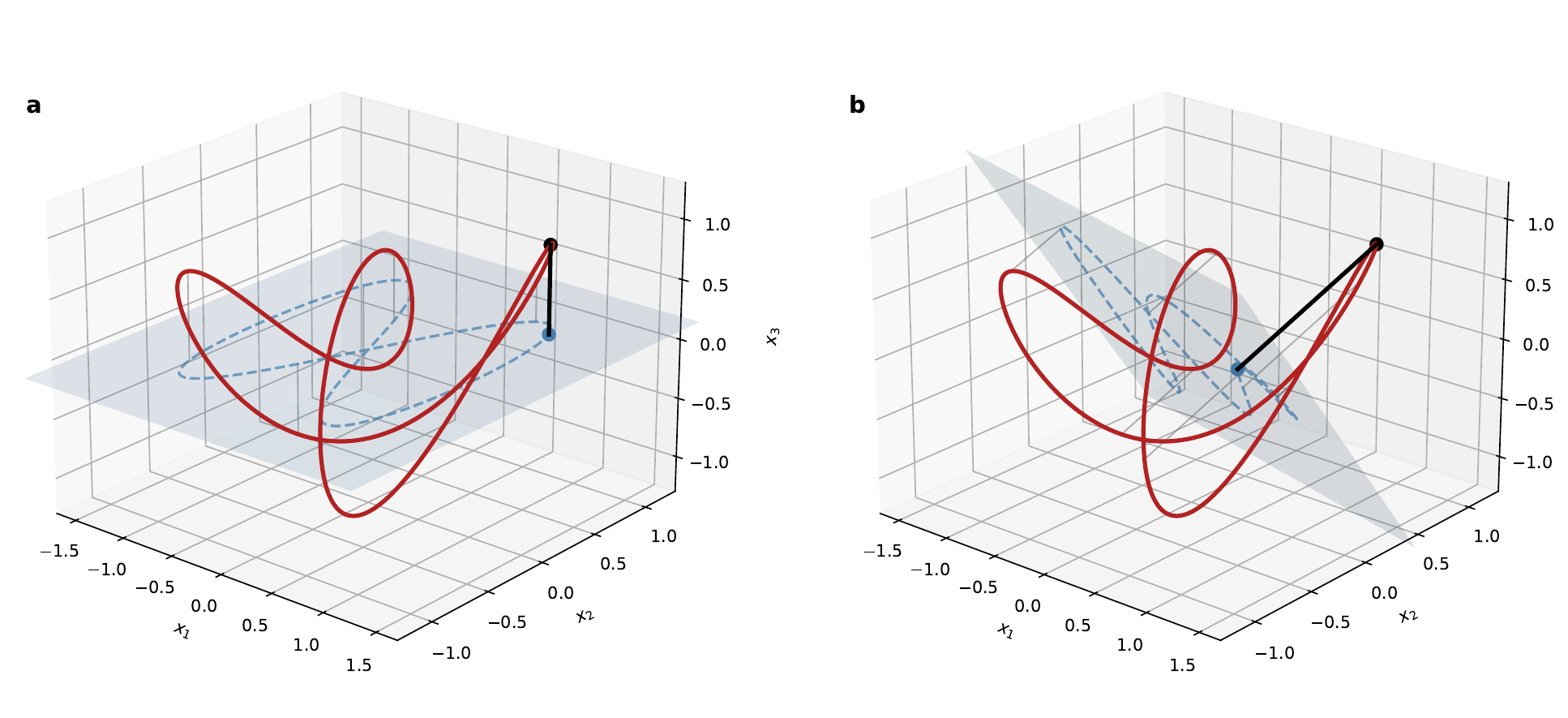}
    \caption{Geometric interpretation of the Kolmogorov $n$-width. The red curve
    is the solution manifold $\mathcal{M}$ embedded in $V = \mathbb{R}^3$, and
    the blue shaded plane is a two-dimensional subspace $V_2$. The dashed blue
    curve is the orthogonal projection of $\mathcal{M}$ onto $V_2$; the  black stick marks the largest perpenducular distance,
    $\sup_{\mu}\|\uu(\mu) - \Pi_{V_2}\uu(\mu)\|$. This is shown for two different planes.}
    \label{fig:nwidth_geometry}
\end{figure}

The practical significance of the $n$-width lies in its role as a universal
lower bound. Suppose one has \emph{any} computational method that produces an
approximation to $\uu(\mu)$ by projection onto an $n$-dimensional subspace
$V_n$. The
worst-case error of this method, by definition, satisfies~\cite{Pinkus1985}:
\begin{equation}
  \label{eq:lower-bound}
  \sup_{\mu \in \mathcal{P}} \norm{\uu(\mu) - \Pi_{V_n}\uu(\mu)}_V
  \geq d_n(\calM, V),
\end{equation}
where $\Pi_{V_n}$ denotes the orthogonal projection onto $V_n$. The left hand side of the inequality above is graphically depicted in \Cref{fig:nwidth_geometry}. No method
based on linear projection onto an $n$-dimensional subspace can beat $d_n$. This gives $d_n$ the character of an information-theoretic limit: it quantifies
the intrinsic compressibility of the parametric problem. If $d_n$ is small for
moderate $n$, the problem is highly amenable to reduced-order modeling, and any
reasonable method should achieve good accuracy with few degrees of freedom. If
$d_n$ remains large even for substantial $n$, no linear reduction method can
be efficient, and one must either accept higher-dimensional approximations or
turn to nonlinear approximation strategies. A method whose
error matches $d_n$ up to a moderate constant factor is \emph{near-optimal}
and cannot be substantially improved by any other linear approach. 

\subsection{Decay rates}
\label{sec:nwidth:decay}

The rate at which $d_n(\calM, V)$ decreases with $n$ determines how
compressible the solution manifold is. Two qualitatively different behaviors
are possible.

\paragraph{Exponential decay.}
For many parametric problems of elliptic type, the $n$-width decays
exponentially:
\begin{equation}
  \label{eq:exp-decay}
  d_n(\calM, V) \leq C\,e^{-\beta\, n^{1/p}},
\end{equation}
where $C$ and $\beta$ are positive constants and $p$ is the dimension of the
parameter space. In this regime, the number of basis functions needed to
achieve an error tolerance $\varepsilon$ grows only as
$n \sim |\log\varepsilon|^p$, which is very favorable. Even for
three-dimensional parameter spaces, one typically needs at most a few dozen
basis vectors for engineering accuracy. The constant $\beta$ controls the
practical efficiency: larger $\beta$ means fewer basis vectors are needed for a
given tolerance. When this exponential behavior holds, it indicates that the solution manifold,
despite living in a high-dimensional ambient space, is effectively a very
low-dimensional object that can be captured by a small subspace. This is the
favorable regime for all reduced-order methods.

This result was first established in the context of single parameter symmetric coercive elliptic PDEs~\cite{Maday2002APC,Maday2002GlobalAP} and subsequently generalized to high-dimensional parameter spaces~\cite{CohenDeVoreSchwab2011,BachmayrCohen2017} where the decay rate was directly linked to the size of the holomorphic extension of the parameter-to-solution map.

\paragraph{Algebraic decay.}
For other problems, notably those involving wave propagation with
discontinuities or transport-dominated phenomena, the $n$-width decays only
algebraically:
\begin{equation}
  \label{eq:alg-decay}
  d_n(\calM, V) \sim C\,n^{-s}
\end{equation}
for some exponent $s > 0$. Such a situation arises, for example, for the linear transport equation
parametrized by velocity~\cite{OhlbergerRave2016} and the second-order wave equation with
discontinuous initial conditions~\cite{GreifUrban2019}.

\paragraph{The regularity connection.}
The dichotomy between exponential and algebraic decay is governed by the
regularity of the parameter-to-solution map $\mu \mapsto \uu(\mu)$. The
fundamental principle is---\emph{the smoother the dependence of the solution on the parameter, the
    faster the $n$-width decays}~\cite{CohenDeVoreSchwab2011,CohenDeVore2015}.
More precisely, if the map $\mu \mapsto \uu(\mu)$ is analytic, then the $n$-width decays exponentially. 


\subsection{Parametric regularity and the role of holomorphy}
\label{sec:nwidth:holomorphy}

The connection between analytic parameter dependence and exponential $n$-width
decay has been made rigorous through a body of work in the reduced basis and
parametric PDE communities.

\paragraph{The size of the holomorphic extension.}
The exponential rate $\beta$ appearing in \eqref{eq:exp-decay} is controlled
by how far the solution map can be analytically continued into the complex
parameter plane~\cite{CohenDeVoreSchwab2011,CohenDeVore2015}. The classical result that makes this precise is Bernstein's
theorem~\cite{bernstein1912ordre}. Consider approximating a function $f: [-1,1] \to \R$ by polynomials
of degree~$n$. If $f$ is analytic on $[-1,1]$, it extends to a holomorphic
function $f(z)$ for complex~$z$ in some neighborhood of the interval. This
extension is valid until one encounters a singularity. The \emph{Bernstein ellipse} $E_\rho$ is
the largest confocal ellipse (with foci at $\pm 1$) whose interior is free of
singularities, and Bernstein's theorem states that the best polynomial
approximation error decays as $\rho^{-n}$, giving $\beta = \log\rho$. 

The
elliptical geometry arises from the Joukowski map $x = (z + z^{-1})/2$ that
relates polynomial approximation on an interval to Fourier analysis on a
circle; for periodic parameter domains---such as the Brillouin zone---the
analogous region is a strip of half-width~$\tau$ in the complex plane, and the
rate is simply $\beta = \tau$, the imaginary distance to the nearest
singularity~\cite{Trefethen2019ApproximationTA}. The interpretation is direct: the farther the nearest complex singularity from
the real parameter domain, the larger~$\beta$, the faster the $n$-width
decays, and the fewer basis vectors are needed. Conversely, as a singularity
approaches the real domain, the rate degrades and more basis vectors are required.

For multivariate parameter domains ($p > 1$), analogous results hold using
multivariate polynomial approximation on polydiscs or polyellipses, leading to
the general bound \eqref{eq:exp-decay} with the characteristic exponent
$n^{1/p}$~\cite{CohenDeVoreSchwab2011,CohenDeVore2015}.

\subsection{Computable bounds via SVD}
\label{sec:nwidth:svd}

The Kolmogorov $n$-width is defined by an infimum over all $n$-dimensional
subspaces and is generally not computable in closed form. However, due to the Eckart-Young theorem\cite{Eckart1936TheAO}, an upper
bound can be obtained from a discrete sampling of the solution manifold using
the singular value decomposition (SVD).

Suppose the solution $\uu(\mu)$ has been computed at a finite collection of
parameter values $\mu_1, \ldots, \mu_M \in \mathcal{P}$, and assemble the
\emph{snapshot matrix}
\begin{equation}
  \label{eq:snapshot-matrix}
  S = \bigl[\uu(\mu_1) \;\; \uu(\mu_2) \;\; \cdots \;\; \uu(\mu_M)\bigr]
  \in \R^{N \times M},
\end{equation}
whose columns are the computed solutions. The SVD of $S$ yields orthonormal
left singular vectors $\boldsymbol{\psi}_1, \boldsymbol{\psi}_2, \ldots$ and
singular values $\sigma_1 \geq \sigma_2 \geq \cdots \geq 0$. By the
Eckart--Young theorem, the subspace
$V_n = \mathrm{span}\{\boldsymbol{\psi}_1, \ldots, \boldsymbol{\psi}_n\}$
minimizes the projection residual over all $n$-dimensional subspaces in the
Frobenius sense. Moreover, the worst-case projection error of
$V_n$ over the snapshot set satisfies
\begin{equation}
  \label{eq:svd-upper}
  \max_{1 \leq \ell \leq M}
  \norm{\uu(\mu_\ell) - \Pi_{V_n}\uu(\mu_\ell)}_V
  \leq \sigma_{n+1},
\end{equation}
since the operator-norm residual of the best rank-$n$ approximation to $S$
equals $\sigma_{n+1}$. Since $V_n$ is one particular $n$-dimensional subspace,
this provides an upper bound on the $n$-width of the discrete snapshot set,
and for sufficiently dense sampling, on the $n$-width of the continuous
manifold $\calM$. When the solution
map $\mu \mapsto \uu(\mu)$ is smooth and the sampling resolves the manifold
adequately, the singular values $\sigma_{n+1}$ provide a tight and
inexpensive proxy for the $n$-width decay rate.

It is worth clarifying the precise relationship between the singular values
of the snapshot matrix and the Kolmogorov $n$-width of the continuous
manifold. Two distinct gaps are involved:

The first gap is between the discrete and continuous $n$-widths. For any
fixed $n$-dimensional subspace $V_n$, the worst-case projection error over
the continuous manifold $\calM$ is at least as large as the worst-case error
over the finite snapshot set $\calM_{\mathrm{disc}} = \{\uu(\mu_1), \ldots,
\uu(\mu_M)\}$, since the supremum is taken over a larger set. Therefore,
\begin{equation}
  \label{eq:discrete-vs-continuous}
  d_n(\calM_{\mathrm{disc}}, V) \leq d_n(\calM, V).
\end{equation}
For a smooth
manifold with sufficiently dense sampling, the gap between the two is
small---specifically,
$d_n(\calM, V) \leq d_n(\calM_{\mathrm{disc}}, V) +
\varepsilon_{\mathrm{samp}}$, where $\varepsilon_{\mathrm{samp}}$ is
controlled by the sampling density and the smoothness of the
parameter-to-solution map.

The second gap is between the discrete $n$-width and the SVD. Even for the
finite snapshot set, the two measure different things: the $n$-width finds
the subspace minimizing the worst-case projection error, whereas the SVD
finds the subspace minimizing the sum of squared projection errors (the
Frobenius residual). The Eckart--Young theorem guarantees that the
worst-case error of the SVD-optimal subspace is bounded by $\sigma_{n+1}$,
so
\begin{equation}
  \label{eq:svd-vs-discrete}
  d_n(\calM_{\mathrm{disc}}, V) \leq \sigma_{n+1},
\end{equation}
but the minimax-optimal subspace could in principle achieve a strictly
smaller worst-case error. 

For sufficiently dense sampling of a smooth manifold, both gaps close and
the singular value sequence $\{\sigma_{n+1}\}$ becomes a tight proxy for
the $n$-width decay rate. \textbf{In practice}, the computation proceeds as follows. One solves the full problem
at a sufficiently dense set of parameter values $\mu_1, \ldots, \mu_M$ and
assembles the snapshot matrix $S$ whose columns are the solutions. The SVD of
$S$ is then computed, yielding the singular values $\sigma_1 \geq \sigma_2
\geq \cdots$. By the bound \eqref{eq:svd-upper}, the sequence
$\{\sigma_{n+1}\}$ provides an upper bound on the $n$-width at each~$n$. Plotting $\log\sigma_n$ versus~$n$ reveals the decay
regime directly. For the bound to be a
reliable proxy for the continuous $n$-width, the sampling must be dense enough
that the discrete snapshot set resolves the solution manifold.

\subsection{Constructive near-optimality: the greedy algorithm}
\label{sec:nwidth:greedy}

The SVD characterization of the preceding subsection requires that the full
problem be solved at every sample point---it is a diagnostic tool, not a
construction method. In practice, one wants to build a good subspace $V_n$
incrementally, using as few full solves as possible.

The greedy algorithm~\cite{Prudhomme2002ReliableRS,veroy2003posteriori,quarteroni2015reduced,Hesthaven2015CertifiedRB} accomplishes this by a sequential worst-case selection.
Starting from an initial subspace $V_1$, the algorithm identifies the
parameter value $\mu^* \in \mathcal{P}$ at which the current subspace
performs worst, solves the full problem at $\mu^*$, and enriches the basis
with the resulting solution. At step $n$, the selection is
\begin{equation}
  \label{eq:greedy-select}
  \mu^* = \arg\max_{\mu \in \mathcal{P}}\;
  \inf_{\mathbf{v} \in V_{n-1}} \norm{\uu(\mu) - \mathbf{v}}_V,
\end{equation}
and the new basis vector is $\uu(\mu^*)$ after orthogonalization against
$V_{n-1}$. In the form \eqref{eq:greedy-select}, the algorithm requires
evaluating the true error at every candidate parameter value, which is as
expensive as the full computation one is trying to avoid. The practical
variant replaces the true error by an inexpensive surrogate---typically a
residual norm that can be evaluated using only the reduced solution---and is
termed the \emph{weak greedy}~\cite{BinevCohenDahmenDeVorePetrovaWojtaszczyk2011}.

Binev et al.~\cite{BinevCohenDahmenDeVorePetrovaWojtaszczyk2011}
proved that the greedy algorithm is \emph{rate-optimal}. The guarantee extends to the weak greedy provided the error surrogate
is \emph{reliable and efficient}~\cite{BinevCohenDahmenDeVorePetrovaWojtaszczyk2011,Hesthaven2015CertifiedRB}. The total
cost is one full solve per greedy step, making the algorithm practical for
problems where each solve is expensive but the number of basis vectors
needed is small. The reduced basis framework has been extended to parametric eigenvalue
problems, where the posteriori error analysis must account for the
spectral gap in the eigenvector error
bound~\cite{Fumagalli2016ReducedBA,Horger2015SimultaneousRB}.

It is worth emphasizing that the greedy algorithm builds the subspace $V_n$
from solution snapshots $\uu(\mu^*)$, but the approximate eigenvalues are
not read off from these snapshots directly. Rather, the original eigenvalue
problem is projected onto $V_n$, and the resulting small system is solved at
each new parameter value. The accuracy of the approximate eigenvalues
depends only on how well $V_n$ contains the true solution---any subspace
that approximates the solution manifold $\calM$ to a given tolerance will
yield eigenvalues of comparable accuracy, regardless of which particular
vectors were used to construct it. This distinction between the basis
vectors and the solutions they represent will become important in the
multi-band analysis of \Cref{sec:analysis:crossings}, where the vectors
used to prove the $n$-width bound differ from the eigenvectors used in
computation, yet both lead to the same subspace and hence the same
eigenvalue accuracy.

\section{$n$-Width Analysis for the Bloch Eigenvalue Problem}
\label{sec:analysis}

We now apply the Kolmogorov $n$-width framework of \Cref{sec:nwidth} to the
Bloch eigenvalue problem of \Cref{sec:bloch:formulation}.  

\paragraph{Holomorphy of the operator family.}
Recall from \eqref{eq:affine-decomp} that the Bloch-transformed stiffness
operator admits the decomposition
$K(\kk) = K_0 + \sum_{m=1}^{Q} f_m(\kk)\,K_m$, where each $f_m(\kk)$ is a
trigonometric polynomial in $\kk$, built from the phase factors
$e^{i\kk\cdot\mathbf{a}_j}$. The same structure holds for the mass operator
$M(\kk)$. Since the exponential function $\kk \mapsto e^{i\kk\cdot\mathbf{a}_j}$
extends to an entire function of $\kk \in \C^d$, the coefficient functions $f_m$ are themselves entire.
It follows that the operator-valued maps $\kk \mapsto K(\kk)$ and
$\kk \mapsto M(\kk)$ extend to holomorphic families of type~(A) in the sense
of Kato~\cite{Kato1995} on all of $\C^d$. In particular, the operators
introduce no singularities: any obstruction to the holomorphy of the eigenpairs
must originate from the eigenvalue problem itself.

\paragraph{Eigenvalue coalescence as a singularity.}
The operators $K(\kk)$ and $M(\kk)$ are entire, so any obstruction to the
holomorphy of the eigenpairs must originate from the eigenvalue problem
itself. The mechanism is the following. Near a point
$\kk_0 \in \C^d$ where two eigenvalues coalesce,
$\omega_i^2(\kk_0) = \omega_j^2(\kk_0)$, the characteristic polynomial
locally has a double root. Generically, this double root unfolds as a
square-root branch point: along any line through $\kk_0$ in the complex
$\kk$-plane, the two eigenvalue branches behave as
$\omega_{i,j}^2(\kk) \approx \lambda_0 \pm c\sqrt{\kk - \kk_0}$, and the
eigenvectors inherit the same branching. No single-valued holomorphic
continuation of either eigenvalue or eigenvector exists around such a point.
On the real parameter domain, Hermitian symmetry provides a special
protection: the eigenvalues of a Hermitian family can always be labeled as
real-analytic functions of the real parameter (Rellich's
theorem~\cite{Kato1995}), and real crossings are not branch points. This
protection is lost once $\kk$ moves into the complex plane, where the
operators are no longer Hermitian~\cite{lu2018level}. Consequently, the nearest complex-$\kk$
coalescence point is a genuine branch point singularity, and it is the
distance from the real parameter domain to these complex branch points that
ultimately controls the holomorphy radius and hence the $n$-width decay rate. The spectral gap on the real domain provides a quantitative lower bound on this distance, as we now make precise.

\paragraph{Holomorphy of eigenpairs via Kato's theory.}
The analyticity of individual eigenpairs is governed by the spectral gap. For a
given band index $j$ and a point $\kk_0 \in \BZ$, define the spectral gap
\begin{equation}
  \label{eq:spectral-gap}
  \delta_j(\kk_0) = \min_{i \neq j}
  \bigl|\omega_j^2(\kk_0) - \omega_i^2(\kk_0)\bigr|,
\end{equation}
which measures the separation of the $j$-th eigenvalue from all other
eigenvalues at the point $\kk_0$. When $\delta_j(\kk_0) > 0$, the eigenvalue $\omega_j^2(\kk_0)$ is simple, and
Kato's analytic perturbation theory~\cite{Kato1995} guarantees that both the
eigenvalue and the eigenvector can be continued as holomorphic functions of
$\kk$ in a neighborhood of $\kk_0$.

As $\kk$ moves away from $\kk_0$ into the complex plane, the eigenvalues
shift at a rate controlled by the variation of the operators. Because
$K(\kk)$ and $M(\kk)$ depend on $\kk$ through the phase factors
$e^{i\kk\cdot\mathbf{a}_j}$ via the affine decomposition
\eqref{eq:affine-decomp}, standard
eigenvalue perturbation theory for Hermitian matrix pencils then gives
\begin{equation}
  \label{eq:eval-lipschitz}
  \bigl|\omega_i^2(\kk) - \omega_i^2(\kk_0)\bigr| \leq L\,\|\kk - \kk_0\|
\end{equation}
for all eigenvalue indices $i$ and all $\kk$ in a neighborhood of $\kk_0$,
where $L$ is a constant which can be computed. For the $j$-th eigenvalue to collide with any other eigenvalue, the total
displacement of the two eigenvalues must bridge the gap $\delta_j(\kk_0)$,
which requires $\|\kk - \kk_0\|$ to be at least of order
$\delta_j(\kk_0)/(2L)$. This gives us the concept of holomorphy radius at $\kk_0$:
 
\begin{proposition}[Holomorphy radius for Bloch eigenpairs]
  \label{prop:holomorphy-radius}
  Let $j$ be a band index and $\kk_0 \in \BZ$ a point at which
  $\omega_j^2(\kk_0)$ is a simple eigenvalue with spectral gap
  $\delta_j(\kk_0) > 0$ as defined in \eqref{eq:spectral-gap}. Let $L > 0$ be
  the Lipschitz constant for the eigenvalue curves, as in
  \eqref{eq:eval-lipschitz}. Then the eigenpair
  $(\omega_j^2(\kk), \uu_j(\kk))$ extends to a holomorphic function of $\kk$
  on the complex ball
  \begin{equation}
    \label{eq:holomorphy-ball}
    B\bigl(\kk_0,\;\rho_j(\kk_0)\bigr)
    = \bigl\{\kk \in \C^d : \|\kk - \kk_0\| < \rho_j(\kk_0)\bigr\},
  \end{equation}
  where the holomorphy radius satisfies the lower bound
  \begin{equation}
    \label{eq:rho-bound}
    \rho_j(\kk_0) \geq \frac{\delta_j(\kk_0)}{2L}.
  \end{equation}
\end{proposition}

The holomorphy radius $\rho_j(\kk_0)$ is a local quantity that varies across
the Brillouin zone. It is large where the $j$-th band is well separated from
its neighbors and small where the spectral gap narrows. The bound \eqref{eq:rho-bound} degrades as $\delta_j \to 0$, and it
breaks down entirely at points of eigenvalue degeneracy, where the
holomorphic continuation of the eigenvector ceases to exist. For the $n$-width
analysis in the following subsection, the relevant quantity is the global
holomorphy radius
\begin{equation}
\label{eq:global-rho}
\rho_j^* = \inf_{\kk \in \BZ}\,\rho_j(\kk)
\geq \frac{1}{2L}\delta_j^*,
\end{equation}
where $\delta_j^* \coloneqq \inf_{\kk \in \BZ}\,\delta_j(\kk)$ is the minimum spectral gap over the entire part of the Brillouin
zone of interest.

\subsection{$n$-Width bounds for isolated bands}
\label{sec:analysis:isolated}

We now combine the holomorphy result of the preceding subsection with the
approximation-theoretic machinery of \Cref{sec:nwidth} to obtain explicit
$n$-width bounds for the solution manifold of an isolated band.

\paragraph{Isolated band.}
We say that the $j$-th band is \emph{isolated} if the spectral gap remains positive over the entire
parameter space of interest: $\delta_j^*>0$. Under this assumption, \Cref{prop:holomorphy-radius} guarantees that the
eigenvector map $\kk \mapsto \uu_j(\kk)$ extends holomorphically to a complex
neighborhood of every point in $\BZ$, with a holomorphy radius that is bounded
below uniformly:
\begin{equation}
  \label{eq:global-rho-isolated}
  \rho_j^* 
  \geq \frac{\delta_j^*}{2L} > 0.
\end{equation}
It follows that $\kk \mapsto \uu_j(\kk)$ extends to a holomorphic function on
the ``thickened'' Brillouin zone
\begin{equation}
  \label{eq:thickened-bz}
  \BZ_{\rho} = \bigl\{\kk \in \C^d : \mathrm{dist}(\kk, \BZ) < \rho\bigr\}
\end{equation}
for any $\rho < \rho_j^*$, and the extension is bounded:
$\sup_{\kk \in \BZ_\rho}\norm{\uu_j(\kk)}_V < \infty$.

\paragraph{Exponential $n$-width decay.}
The solution manifold for the $j$-th band,
$\calM_j = \{\uu_j(\kk) \in V : \kk \in \BZ\}$, is the image of the
Brillouin zone under the holomorphic map $\kk \mapsto \uu_j(\kk)$.
Since this map is periodic on $\BZ$ and extends holomorphically to the
strip $\BZ_\rho = \{\kk \in \C^d : \mathrm{dist}(\kk, \BZ) < \rho\}$
for any $\rho < \rho_j^*$, multivariate Fourier approximation estimates
on periodic domains apply directly: the best $n$-term trigonometric
approximation error decays exponentially with a rate controlled by the
strip half-width~$\rho_j^*$. This is the periodic-domain analog of the
polynomial approximation results established by Cohen and
DeVore~\cite{CohenDeVore2015} and Bachmayr and
Cohen~\cite{BachmayrCohen2017} for parametric elliptic source problems;
the same approximation-theoretic principles carry over to the present
eigenvalue setting because Kato's theory furnishes the required
holomorphic parameter-to-solution map. Combining these estimates with the
$n$-width machinery of \Cref{sec:nwidth:holomorphy}, we obtain the
following result.

\begin{theorem}[$n$-width bound for an isolated band]
  \label{thm:isolated-band}
  Let $j$ be the index of an isolated band, and let $\calM_j$ be the corresponding solution
  manifold \eqref{eq:sol-manifold}. Then the Kolmogorov $n$-width of $\calM_j$
  in $V$ satisfies
  \begin{equation}
    \label{eq:nwidth-isolated}
    d_n(\calM_j, V) \leq C\,e^{-\beta\,n^{1/d}},
  \end{equation}
  where $d = \dim(\kk)$ is the dimension of the Brillouin zone and the
  constants are given by:
  \begin{itemize}[nosep]
    \item $\beta > 0$ depends on the global holomorphy radius $\rho_j^*$ and
      hence, via \eqref{eq:global-rho-isolated}, on 
      $\delta_j^*,L$---$\beta$ grows as $\delta_j^*$ increases or $L$ decreases;
    \item $C > 0$ is directly proportional to
      $\sup_{\kk \in \BZ_\rho}\norm{\uu_j(\kk)}_V$, the norm of the
      holomorphic extension.
  \end{itemize}
\end{theorem}

\Cref{thm:isolated-band} provides a quantitative justification for the
empirical success of reduced-order methods applied to well-separated bands. The
exponential bound \eqref{eq:nwidth-isolated} implies that the number of basis
vectors needed to approximate the solution manifold $\calM_j$ to within a
tolerance $\varepsilon$ scales as
\begin{equation}
  \label{eq:basis-count}
  n(\varepsilon) \sim \Bigl(\frac{1}{\beta}\log\frac{C}{\varepsilon}\Bigr)^d.
\end{equation}
For a two-dimensional crystal ($d = 2$) with a well-separated band, achieving
engineering accuracy ($\varepsilon \sim 10^{-4}$) may require only
$n \sim |\log\varepsilon|^2 \approx 100$ or fewer basis vectors, compared to a
full discretization dimension $N$ in the tens of thousands. The bound also
reveals that bands
that are well separated from their neighbors are more compressible, while bands
with narrow gaps require more basis vectors for the same accuracy.

\subsection{Band crossings, spectral projectors, and multi-band formulations}
\label{sec:analysis:crossings}

The exponential $n$-width bound of \Cref{thm:isolated-band} rests on the
isolated band assumption. In practice, band structures
almost always exhibit crossings or near-crossings somewhere in the Brillouin
zone. In this subsection we analyze how band crossings affect the $n$-width
decay, and show that working with spectral projectors onto groups of bands
rather than individual eigenvectors resolves all forms of degeneracy in a
unified way.

\subsubsection{Avoided crossings}
\label{sec:analysis:crossings:avoided}

An avoided crossing occurs when two bands approach each other closely but do
not touch - $\delta_j^*$ remains positive but can be arbitrarily small. This is the generic
behavior for bands not protected by symmetry. The isolated band assumption
still holds, so \Cref{thm:isolated-band} continues to
apply; the exponential bound \eqref{eq:nwidth-isolated} is not invalidated.
However, the rate $\beta$ degrades because 
$\rho_j^*$ shrinks in proportion to $\delta_j^*$. With
$\beta = \mathcal{O}(\delta_j^*/L)$, the required basis dimension scales as
\begin{equation}
  \label{eq:avoided-basis-count}
  n(\varepsilon) \sim \left(\frac{L}{\delta_j^*}
  \log\frac{C}{\varepsilon}\right)^d,
\end{equation}
which diverges as $\delta_j^* \to 0$. The decay remains exponential but with
an arbitrarily poor rate as the crossing tightens.

\subsubsection{Exact degeneracies and the spectral projector}
\label{sec:analysis:crossings:projectors}

At an exact band crossing, individual eigenvectors $\uu_j(\kk)$ cease to be
the right objects entirely. Two problems arise simultaneously. First, there is
a \emph{gauge ambiguity}: the eigenvector is defined only up to a scalar
multiple. Second, and more fundamentally, at a degeneracy $\uu_j(\kk)$ is not
even well-defined as a continuous function of $\kk$: approaching the crossing
from different directions, it may converge to different elements of the
degenerate eigenspace. The holomorphic continuation of
\Cref{prop:holomorphy-radius} breaks down.

Both difficulties are resolved by replacing individual eigenvectors with the
\emph{spectral projector} onto the group of $J$ bands of interest:
\begin{equation}
  \label{eq:multiband-proj}
  P_J(\kk) = \frac{1}{2\pi i} \oint_{\Gamma_J}
  \bigl(K(\kk) - \zeta\,M(\kk)\bigr)^{-1} M(\kk)\,\dd\zeta,
\end{equation}
where $\Gamma_J \subset \C$ encloses $\omega_1^2(\kk), \ldots,
\omega_J^2(\kk)$ but no other eigenvalues. The operator $P_J(\kk)$ projects
orthogonally onto the $J$-dimensional spectral subspace
\begin{equation}
  \label{eq:spectral-subspace}
  \calS_J(\kk) = \mathrm{range}\bigl(P_J(\kk)\bigr)
  = \mathrm{span}\bigl\{\uu_1(\kk), \ldots, \uu_J(\kk)\bigr\}.
\end{equation}
The spectral projector is gauge-invariant by construction, and insensitive to
crossings within the group: rearrangements of eigenvalues inside $\Gamma_J$
leave the contour integral \eqref{eq:multiband-proj} unchanged, so
$\calS_J(\kk)$ varies holomorphically through any internal crossing. The only
event that disrupts the projector is an eigenvalue crossing at the cluster
boundary---the $J$-th meeting the $(J+1)$-th.

\paragraph{Symmetry-enforced degeneracies and conical intersections.}
Both types of exact degeneracy arising in phononic band structures are handled
uniformly by this framework. At symmetry-enforced degeneracies, point-group
symmetry forces $m \geq 2$ eigenvalues to coalesce exactly; individual
eigenvectors are ill-defined, but $P_J(\kk)$ remains holomorphic provided the
cluster is separated from the remaining spectrum. At conical intersections
(Dirac points), two bands touch with linear dispersion 
and the eigenvectors acquire a Berry phase $\phi_B = \pi$ upon encircling the
intersection~\cite{Berry1984QuantalPF}. This sign reversal is a property of the
\emph{eigenvector frame} within the two-dimensional subspace; it does not
obstruct the smoothness of the \emph{subspace itself}. The spectral projector, measuring the subspace rather than the
frame, therefore remains holomorphic through a conical intersection by the same
contour argument. 

In both cases the conclusion is the same: as long as the cluster of $J$ bands
is separated from the $(J+1)$-th band by a positive gap, $P_J(\kk)$ is
holomorphic on a complex neighborhood of $\BZ$, and the $n$-width analysis
proceeds exactly as for an isolated band, with $P_J(\kk)$ in place of
$\uu_j(\kk)$ and the cluster gap in place of the individual band gap.

\subsubsection{Multi-band $n$-width bound}
\label{sec:analysis:multiband}

To apply the $n$-width framework of \Cref{def:nwidth}, we need a compact
subset of~$V$, not an operator-valued map. When the gap between the $J$-th
and $(J+1)$-th bands is uniformly positive over~$\BZ$, Kato's analytic
perturbation theory~\cite{Kato1995} guarantees the existence of an
orthonormal set $\{\uu_1(\kk), \ldots, \uu_J(\kk)\}$ spanning
$\calS_J(\kk)$ that depends holomorphically on~$\kk$ throughout the
thickened Brillouin zone $\BZ_\rho$ for $\rho \lesssim \delta_J^*/(2L)$.
The multi-band solution manifold is:
\begin{equation}
  \label{eq:multiband-sol-manifold}
  \calM_J = \bigcup_{j=1}^{J}
  \bigl\{\uu_j(\kk) \in V : \kk \in \BZ\bigr\}
\end{equation}
with the
cluster gap $\delta_J^*$ in place of the individual band gap. Since any
$n$-dimensional subspace approximating $\calM_J$ must devote at least $J$
dimensions to spanning the $J$-dimensional eigenspace $\calS_J(\kk)$ at any
single $\kk$-point, the remaining $n - J$ dimensions are what capture the
variation across~$\BZ$, giving the following.
 
\begin{theorem}[$n$-width bound for a multi-band manifold]
  \label{thm:multiband}
  Let $J \geq 1$ and suppose the gap between the $J$-th and $(J+1)$-th bands
  satisfies $\delta_J^* = \inf_{\kk \in \BZ}|\omega_J^2(\kk) -
  \omega_{J+1}^2(\kk)| > 0$. Then
  \begin{equation}
    \label{eq:nwidth-multiband}
    d_n(\calM_J, V) \leq C\,e^{-\beta\,(n-J)^{1/d}}, \qquad n \geq J,
  \end{equation}
  where $\beta > 0$ depends on $\delta_J^*/(2L)$ and $C > 0$ depends on
  $\sup_{\kk \in \BZ_\rho}\max_j\norm{\uu_j(\kk)}_V$ and $J$.
\end{theorem}

The central message is that exponential decay holds regardless of what
crossings occur among the first $J$ bands. The relevant quantity is not the
type of degeneracy within the cluster, but simply whether the cluster is
isolated from the bands outside it. All reduction methods---whether based on
sampling eigenvectors at selected $\kk$-points, component mode synthesis, or
snapshot compression---produce a specific $n$-dimensional subspace $V_n$
approximating $\calM_J$, and $d_n(\calM_J, V)$ quantifies the gap between any
such method and the theoretical optimum. Methods that retain all eigenvectors
within a frequency window implicitly work with the full spectral subspace and
should exhibit exponential convergence whenever the window boundaries lie in
spectral gaps.

\subsection{Domain, bands, and resolution}
Every reduced-order band structure computation involves some choices. Together they
determine the solution manifold whose $n$-width is the object of study.

The Brillouin zone $\BZ$ is the natural parameter domain, but one rarely needs
the full zone. In practice, $\mathcal{B}$ might be a high-symmetry path
$\Gamma$--$X$--$M$--$\Gamma$, a subdomain of interest for a particular
application, or the full zone. This choice matters because the exponential
rate $\beta$ is controlled by the minimum spectral gap \emph{over
$\mathcal{B}$}. 

One must decide how many bands to include in the reduced model. This
determines the cluster $\{\omega_1^2(\kk), \ldots, \omega_J^2(\kk)\}$ and, where the cluster boundary falls. The key observation, established
in \Cref{sec:analysis:crossings}, is that \emph{crossings within the cluster
are free}. Only the gap $\delta_J^*$ between the
$J$-th and $(J+1)$-th bands over $\mathcal{B}$ matters for the $n$-width.

An important distinction is between approximating the \emph{spectral
subspace} $\calS_J(\kk)$ spanned by the $J$ bands and tracking
\emph{individual eigenvectors} $\uu_j(\kk)$ within that subspace. For the
large majority of engineering quantities of interest---dispersion surfaces,
group velocities, transmission spectra, density of states---the spectral
subspace is sufficient. Individual band labels are not needed, and the reduced
model need only reproduce the subspace accurately.

When subspace approximation suffices, exact degeneracies within the cluster
are completely harmless: the spectral projector does not distinguish between
crossing and non-crossing bands, and exponential convergence holds as long as
the cluster gap $\delta_J^*$ is positive. The situation in which exponential
convergence is genuinely obstructed is narrow and specific: it requires an
exact degeneracy at some $\kk_0 \in \mathcal{B}$ \emph{and} a requirement to
resolve the individual bands involved. This arises, for instance, when
computing topological band invariants such as the Berry phase or Chern
numbers, which are properties of individual eigenvector bundles rather than of
the subspace. This distinction between subspace reducibility and individual eigenvector
reducibility has been studied extensively in the context of Wannier
functions in electronic structure theory~\cite{MarzariVanderbilt1997,MarzariMostofiYatesSouzaVanderbilt2012}.
The construction of maximally localized Wannier functions---a real-space
basis for the spectral subspace---succeeds precisely when the composite
band group admits a smooth frame, which fails if the Chern number is
nonzero. The spectral projector onto the composite bands remains smooth
regardless, paralleling the distinction made here between the holomorphy of
$P_J(\kk)$ and the potential non-smoothness of individual eigenvector
branches.

Thus, exponential convergence of the $n$-width is the generic
situation for practically relevant computations. The three choices above
determine the effective gap $\delta^*$ that controls the exponential rate
$\beta \sim \delta^*/(2L)$.

\section{Explorations for a 1-D Problem}
\label{sec:1d}

The preceding sections developed the $n$-width framework for phononic band
structures in generality. In this section, we examine the theory's predictions
through a concrete one-dimensional example that is simple enough for exact
computation but rich enough to exhibit the phenomena of interest. Working in
one dimension has two specific advantages: the Brillouin zone is an interval
($d = 1$), so the theory predicts pure exponential $n$-width decay
$d_n \leq C\,e^{-\beta n}$ rather than the stretched exponential
$e^{-\beta n^{1/d}}$ that arises for $d \geq 2$; and all eigenvalues are
generically simple, so the isolated band assumption
(\Cref{sec:analysis:isolated}) holds throughout the Brillouin zone, and the
complications associated with band crossings
(\Cref{sec:analysis:crossings}) do not arise.

\subsection{Problem setup and solution method}
\label{sec:1d:setup}

Consider a one-dimensional phononic crystal consisting of an infinite rod with
periodically varying stiffness $E(x)$ and density $\rho(x)$, both of period
$a$. The time-harmonic scalar wave equation, after application of Bloch's
theorem with wave number $k \in \BZ = [0, \pi/a]$, reduces to an eigenvalue
problem for the periodic part $\tilde{u}(x; k)$ of the Bloch mode on the unit
cell $[0, a]$:
\begin{equation}
  \label{eq:1d-bloch}
  -\bigl(\tfrac{d}{dx} + ik\bigr)
  \Bigl[E(x)\bigl(\tfrac{d}{dx} + ik\bigr)\tilde{u}\Bigr]
  = \omega^2\,\rho(x)\,\tilde{u},
  \qquad \tilde{u}(0) = \tilde{u}(a),
  \quad E(0)\tilde{u}'(0) = E(a)\tilde{u}'(a).
\end{equation}
For each $k$, this problem has a countable sequence of eigenvalues
$0 \leq \omega_1^2(k) \leq \omega_2^2(k) \leq \cdots$, and the
corresponding eigenfunctions $\tilde{u}_j(x; k)$ define the solution manifold
$\calM_j = \{\tilde{u}_j(\cdot\,; k) : k \in \BZ\}$ for the $j$-th band. We take a two-harmonic property profile as the primary example:
\begin{align}
  \label{eq:two-harmonic-profile}
  E(x) &= 1 + \alpha_1\cos\bigl(2\pi x / a\bigr)
           + \alpha_2\cos\bigl(4\pi x / a\bigr), \notag\\
  \rho(x) &= 1 + \alpha_1\cos\bigl(2\pi x / a\bigr)
              - \alpha_2\cos\bigl(4\pi x / a\bigr),
\end{align}
with $\alpha_1 = 0.6$, $\alpha_2 = 0.3$, and $a = 1$(\Cref{fig:phononic_1d}(a)). The eigenvalue problem \eqref{eq:1d-bloch} is solved using the
transfer matrix method (TMM). The unit cell is divided into $N_\ell = 100$
thin sublayers, each treated as homogeneous with properties evaluated at its
midpoint. The $2 \times 2$ transfer matrix of the full cell is formed as the
product of the sublayer transfer matrices, and the dispersion relation
$\cos(ka) = \tfrac{1}{2}\mathrm{tr}\,T(\omega)$ is solved for the
eigenfrequencies $\omega_j(k)$ by root-finding. The eigenfunctions are
obtained by propagating the Bloch eigenvector of $T(\omega_j)$ through the
sublayers and removing the plane-wave factor. Each eigenfunction is normalized
to unit $L^2$ norm on $[0, a]$, and a gauge convention is imposed by
requiring $\tilde{u}_j(0; k)$ to be real and positive. The spatial
discretization on $N_\ell + 1 = 101$ grid points serves as the ambient space
$V = \C^{101}$ for the $n$-width analysis.

\subsection{Snapshot SVD and $n$-width characterization}
\label{sec:1d:svd}

For each band index $j$, we sample the solution manifold $\calM_j$ by
evaluating the eigenfunction $\tilde{u}_j(\cdot\,; k)$ at $M = 200$
uniformly spaced points $k_1, \ldots, k_M$ in the interior of $\BZ$, and
assemble the snapshot matrix
\begin{equation}
  \label{eq:1d-snapshot}
  S_j = \bigl[\tilde{u}_j(\cdot\,; k_1) \;\; \cdots \;\;
  \tilde{u}_j(\cdot\,; k_M)\bigr] \in \C^{101 \times 200},
\end{equation}
with columns weighted by $\sqrt{\Delta x}$ so that the discrete $\ell^2$ norm
of each column approximates the $L^2$ norm of the corresponding
eigenfunction. As discussed in \Cref{sec:nwidth:svd}, the singular values
$\sigma_1^{(j)} \geq \sigma_2^{(j)} \geq \cdots$ of $S_j$ provide an upper
bound on the Kolmogorov $n$-width of $\calM_j$:
\begin{equation}
  \label{eq:1d-svd-bound}
  d_n(\calM_j, V) \leq \sigma_{n+1}^{(j)} + \varepsilon_{\mathrm{samp}},
\end{equation}
with $\varepsilon_{\mathrm{samp}}$ small for the sampling density used.

\begin{figure}[htbp]
    \centering
    \includegraphics[width=\textwidth]{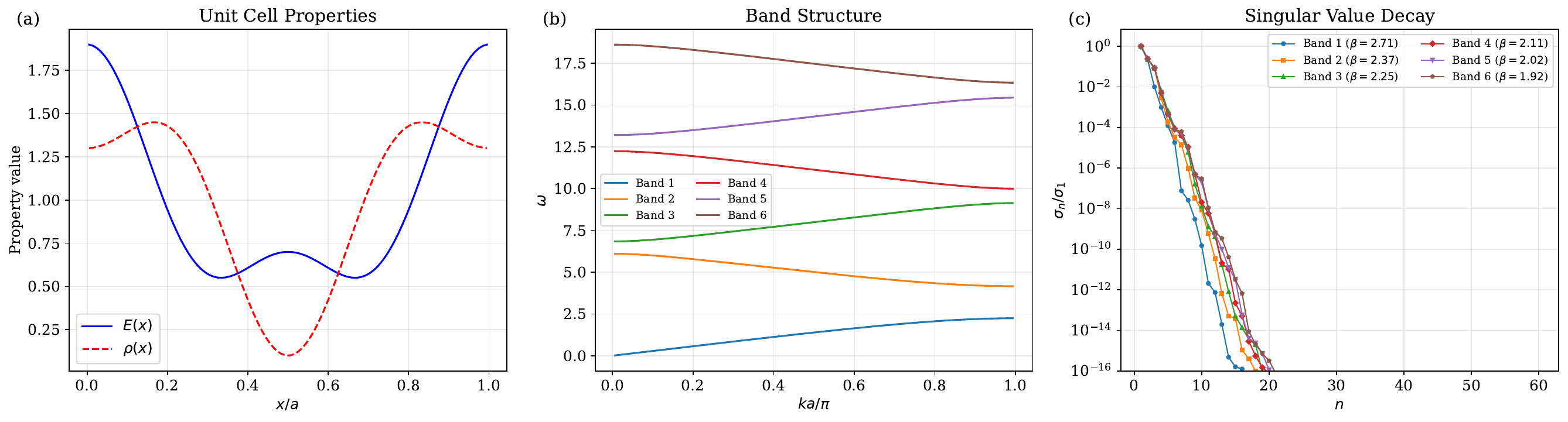}
    \caption{1D phononic crystal with continuously varying properties. (a) Unit cell property distributions $E(x)$ and $\rho(x)$. (b) Band structure showing the first six dispersion branches. (c) Normalized singular value decay of the snapshot matrices for each band, with fitted exponential rates $\beta$.}
    \label{fig:phononic_1d}
\end{figure}

\Cref{fig:phononic_1d}(a) shows the modulus and density profile over the unit cell under consideration. \Cref{fig:phononic_1d}(b) shows the computed phononic bandstructure of this composite for the first six bands. Finally, \Cref{fig:phononic_1d}(c) shows the normalized singular values
$\sigma_n^{(j)}/\sigma_1^{(j)}$ as a function of $n$ for the first six bands.
In every case, the singular values decay exponentially over many orders of
magnitude before reaching machine precision near $n = 15$. The decay is well
described by the model $\sigma_n \sim C\,e^{-\beta n}$ and the legend also includes the estimated $\beta$ for each band.

\subsection{Oracle greedy}
\label{sec:1d:oracle_greedy}

For the greedy experiments of this and the following subsection, we switch
to the single-harmonic profile
\begin{equation}
  \label{eq:sinusoidal-profile}
  E(x) = \rho(x) = 1 + \alpha\cos\bigl(2\pi x / a\bigr),
\end{equation}
with $\alpha = 0.8$, giving a contrast ratio $E_{\max}/E_{\min} = 9$.
The proportionality $E \propto \rho$ simplifies the transfer matrix
structure while preserving the essential features of the $n$-width
analysis. All other parameters remain unchanged.

We now compare the SVD benchmark against the greedy algorithm of
\Cref{sec:nwidth:greedy}, implemented here in its oracle form: at each step,
the true projection error is evaluated over all precomputed snapshots, and the
worst-approximated snapshot is added to the basis. This is the idealized
version of the greedy selection \eqref{eq:greedy-select}---it requires all
$J \cdot M = 2000$ snapshots to be available, so it is not a practical
algorithm, but it provides a clean test of the rate-optimality guarantee
without the complication of an approximate error surrogate.

The basis is initialized with the $J = 10$ eigenvectors at the Brillouin zone
center $k = \pi/(2a)$, providing the minimum $J$-dimensional subspace needed
to represent the spectral subspace $\calS_J(k)$ at a single $k$-point. From
this starting point, the greedy adds one vector per step, selecting from the
full snapshot set.

\begin{figure}[htbp]
    \centering
    \includegraphics[width=\textwidth]{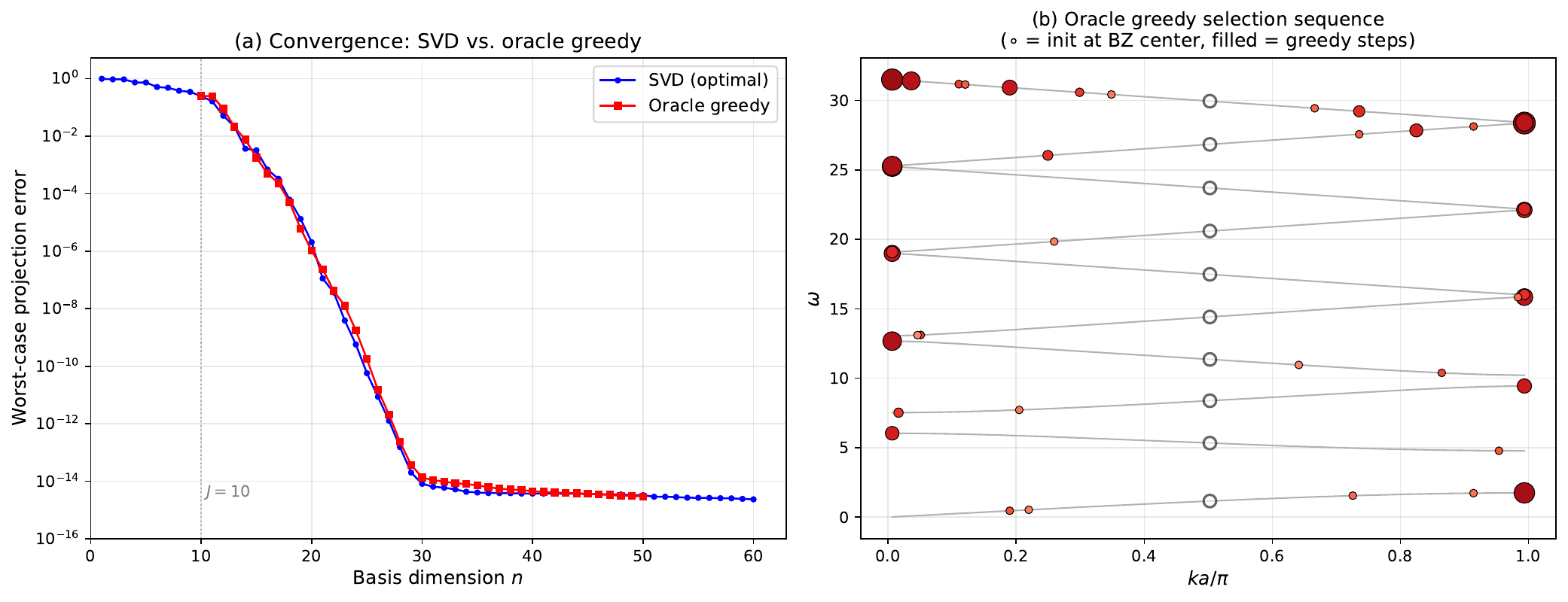}
    \caption{Oracle greedy algorithm for the first $J = 10$ bands.
    \textbf{(a)}~Worst-case projection error versus basis dimension~$n$ for
    the SVD-optimal subspace and the oracle greedy.
    \textbf{(b)}~Band structure with greedy selections marked. Open circles
    denote the initialization point; filled circles indicate subsequent
    greedy selections, with marker size decreasing in the order of selection.}
    \label{fig:oracle}
\end{figure}

\Cref{fig:oracle}(a) shows the worst-case projection error as a function of
the basis dimension~$n$ for both the SVD-optimal subspace and the oracle
greedy. For $n < J = 10$, the SVD error remains large---the subspace cannot
yet span the ten-dimensional eigenspace at any single $k$-point. Beyond
$n = J$, both curves decay exponentially in $n - J$, consistent with the
bound \eqref{eq:nwidth-multiband} of \Cref{thm:multiband}. The greedy tracks
the SVD, confirming the
rate-optimality guarantee of Binev et al.~\cite{BinevCohenDahmenDeVorePetrovaWojtaszczyk2011}: the greedy
achieves the same exponential decay rate as the optimal subspace despite
constructing the basis sequentially. Both methods reach machine precision near
$n = 30$.

\Cref{fig:oracle}(b) reveals the order in which the greedy selects its basis
vectors, providing physical insight into the structure of the solution
manifold. The first post-initialization steps select the highest bands (bands
$8$--$10$) at the zone boundaries $k \approx 0$ and $k \approx \pi/a$. This
is natural: the basis was initialized at the zone center, and the
eigenfunctions at the zone boundaries---which correspond to standing waves of
maximally different character---are the most poorly represented by the initial
subspace. Among the bands, the highest ones are selected first because they generally
have the smallest spectral gaps and hence the fastest-varying eigenvectors,
consistent with the per-band $n$-width analysis of \Cref{sec:1d:svd}. The
lower bands, with their larger gaps and slower variation, converge with fewer
additional vectors and are selected last.

This selection pattern has implications for the design of reduction methods. It
provides a principled justification for the empirical practice of building
reduced bases from eigenvectors computed at high-symmetry points, as is done in the Hussein's Reduced Bloch Mode Expansion (RBME) method~\cite{Hussein2009}: the greedy algorithm independently discovers that these are the
most informative sampling locations. At the same time, it shows that the
zone-boundary strategy becomes suboptimal as the basis grows---the greedy
transitions to selecting interior $k$-points for the higher bands once the
zone-boundary information has been absorbed.

\subsection{Residual-based greedy algorithm}
\label{sec:1d:residual-greedy}

\begin{figure}[htbp]
    \centering
    \includegraphics[width=\textwidth]{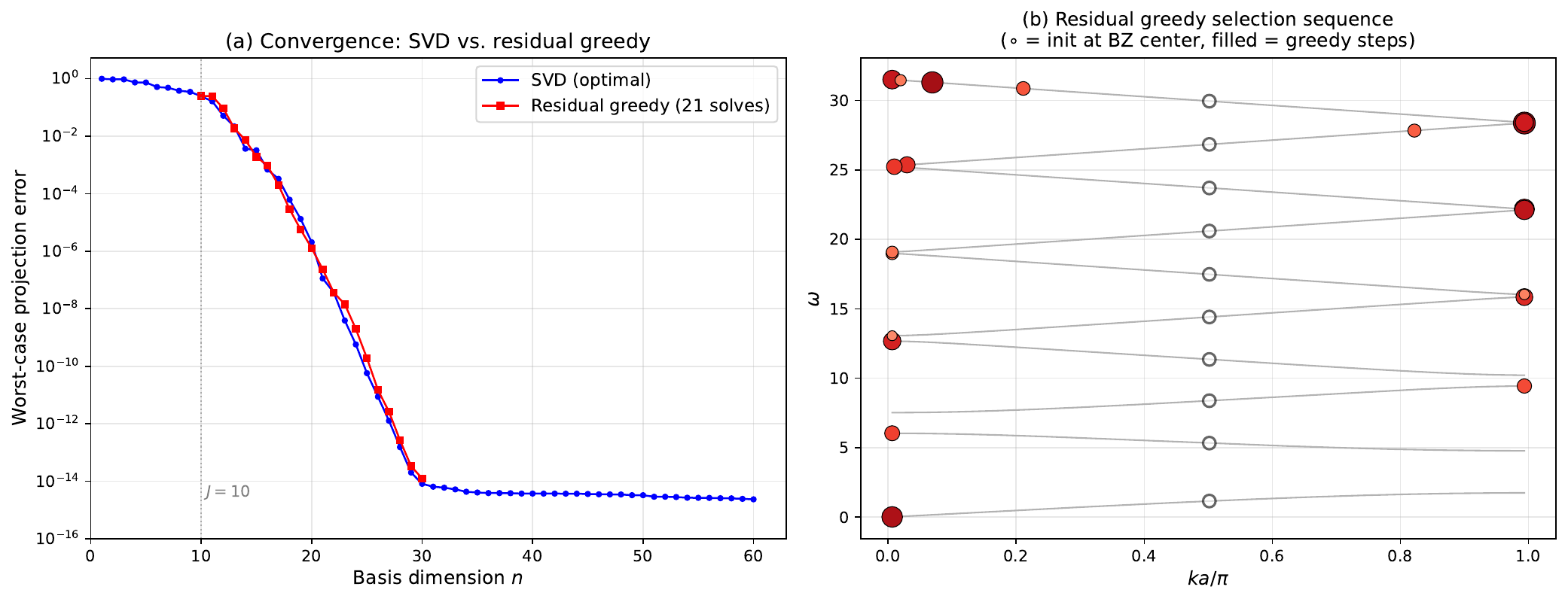}
    \caption{Residual-based greedy algorithm for the first $J = 10$ bands,
    initialized with $J$ eigenvectors at the Brillouin zone center.
    \textbf{(a)}~Worst-case projection error versus basis dimension~$n$ for
    the SVD-optimal subspace and the residual-based greedy. The algorithm
    converges to machine precision at $n = 30$ using $21$ full eigenvalue
    solves.
    \textbf{(b)}~Band structure with greedy selections marked. Open circles
    denote the initialization point; filled circles indicate subsequent
    greedy selections, with marker size decreasing in the order of selection.}
    \label{fig:residual-greedy}
\end{figure}

The oracle greedy of the preceding subsection requires all snapshots to be
computed in advance---it is a diagnostic tool for characterizing the $n$-width,
not a practical algorithm. We now implement the weak greedy of
\Cref{sec:nwidth:greedy} using a residual-based error surrogate that avoids
this requirement, with the goal of verifying that the practical algorithm
achieves the same convergence rate.

The algorithm targets the first $J$ bands simultaneously and constructs the
basis one vector at a time. It is initialized identically to the oracle
greedy: one full eigenvalue solve at the Brillouin zone center, with the
first $J$ eigenvectors forming the starting basis
$\Phi \in \C^{N \times J}$.

At each subsequent step, the algorithm must identify the wave vector $\kk^*$
at which the current basis performs worst. For a general parametric problem,
this would require solving the full problem at every candidate $\kk$-point,
defeating the purpose. The affine decomposition
\eqref{eq:affine-decomp} provides a way around this. Since
$K(\kk) = K_0 + \sum_{m} f_m(\kk)\,K_m$ and likewise for $M(\kk)$, the
operator--basis products $K_m\,\Phi$ and $M_m\,\Phi$ can be precomputed
once per greedy step. The reduced
matrices $\widetilde{K}(\kk) = \Phi^H K(\kk)\,\Phi$ and
$\widetilde{M}(\kk) = \Phi^H M(\kk)\,\Phi$ are then assembled for each
candidate $\kk$ from these precomputed products. The $n \times n$ reduced eigenvalue
problem is solved, yielding approximate eigenpairs
$(\widetilde{\omega}_j^2,\,\widetilde{\uu}_j)$ for $j = 1, \ldots, J$.
The residual of each reduced eigenpair,
\begin{equation}
  \label{eq:residual-def}
  \mathbf{r}_j(\kk)
  = K(\kk)\,\Phi\,\widetilde{\uu}_j
  - \widetilde{\omega}_j^2\,M(\kk)\,\Phi\,\widetilde{\uu}_j,
\end{equation}
is an $N$-vector that measures how well the reduced eigenpair satisfies the
full eigenvalue equation. The error indicator at each
wave vector is the worst residual over all $J$ bands:
\begin{equation}
  \label{eq:greedy-indicator}
  \Delta_n(\kk) = \max_{j=1,\ldots,J} \norm{\mathbf{r}_j(\kk)}.
\end{equation}
The greedy selects $\kk^* = \arg\max_{\kk \in \mathcal{K}}\,\Delta_n(\kk)$
over a dense training grid $\mathcal{K} \subset \BZ$, performs one full
$N$-dimensional eigenvalue solve at $\kk^*$, identifies which of the $J$
eigenvectors at $\kk^*$ is worst-approximated by the current basis,
orthogonalizes it against $\Phi$, and appends it. The basis dimension
increases by one and the process repeats.

\Cref{fig:residual-greedy}(a) shows the convergence of the residual-based
greedy alongside the SVD benchmark. The residual-based algorithm tracks the
SVD closely. Both methods exhibit exponential decay in
$n - J$, consistent with the bound \eqref{eq:nwidth-multiband} of
\Cref{thm:multiband}. The residual-based greedy reaches machine precision at
$n = 30$ using only $21$ full eigenvalue solves.

\Cref{fig:residual-greedy}(b) shows the selection sequence on the band
structure. The qualitative pattern matches the oracle greedy of
\Cref{fig:oracle}: early steps target the highest bands at the zone
boundaries, where the eigenvectors differ most from the zone-center
initialization, and later steps fill in the lower bands. The specific ordering
differs slightly---the residual estimator weights contributions differently
from the true projection error---but the overall convergence rate and the
final basis dimension are comparable.


\section{2-D Problem}
\label{sec:2d}

We now extend the numerical investigation to two dimensions, where the
Brillouin zone is a two-dimensional domain and the theory predicts the
stretched exponential decay $d_n \leq C\,e^{-\beta\,n^{1/2}}$ rather than
the pure exponential that arises in one dimension. The primary objective is
to test whether the exponent $n^{1/d}$ in the $n$-width bound of
\Cref{thm:isolated-band} correctly captures the role of the parameter-space
dimension.

\subsection{Problem setup}
\label{sec:2d:setup}

Consider a two-dimensional phononic crystal with a square lattice of side
$a$. The unit cell $\Omega = [0,a]^2$ contains a circular inclusion of
radius $r$ centered at $(a/2,\, a/2)$, embedded in a matrix material. We
take $a = 1$ and $r/a = 0.35$. The governing equation is the scalar wave equation
\begin{equation}
  \label{eq:2d-scalar}
  -\nabla \cdot \bigl[E(\xx)\,\nabla \tilde{u}\bigr]
  = \omega^2\,\rho(\xx)\,\tilde{u},
\end{equation}
for the Bloch-periodic part $\tilde{u}(\xx;\kk)$ on the unit cell, with
$\tilde{u}$ periodic on $\partial\Omega$. The material properties are
piecewise constant: $E = 1$, $\rho = 1$ in the matrix and $E = 12$,
$\rho = 1$ in the inclusion.

The unit cell is discretized with linear triangular finite elements on a
mesh generated by Gmsh~\cite{geuzaine2009gmsh}, with periodic meshing
constraints imposed to ensure matching nodes on opposite boundaries.
\Cref{fig:2d-mesh-bands}(a) shows the resulting mesh. The global stiffness
and mass matrices $K_{\mathrm{g}}$ and $M_{\mathrm{g}}$ are assembled
without any Bloch boundary conditions. At each wave vector $\kk$, the
Bloch-periodic boundary conditions are imposed through the transformation
matrix $T(\kk)$, which identifies degrees of freedom on opposite edges and
corners of the unit cell through the phase factors
$z_1 = e^{i k_x a}$ and $z_2 = e^{i k_y a}$:
\begin{equation}
  \label{eq:2d-bloch-bc}
  u_{\mathrm{right}} = z_1\,u_{\mathrm{left}}, \qquad
  u_{\mathrm{top}} = z_2\,u_{\mathrm{bottom}}, \qquad
  u_{\mathrm{TR}} = z_1 z_2\,u_{\mathrm{BL}},
\end{equation}
with analogous relations for the remaining corners.
The reduced eigenvalue problem
$\hat{K}(\kk)\,\hat{\uu} = \omega^2\,\hat{M}(\kk)\,\hat{\uu}$, where
$\hat{K} = T^H K_{\mathrm{g}}\,T$ and
$\hat{M} = T^H M_{\mathrm{g}}\,T$, is then solved at each $\kk$-point.

\subsection{Band structure}
\label{sec:2d:bands}

\begin{figure}[htbp]
  \centering
  \includegraphics[width=\textwidth]{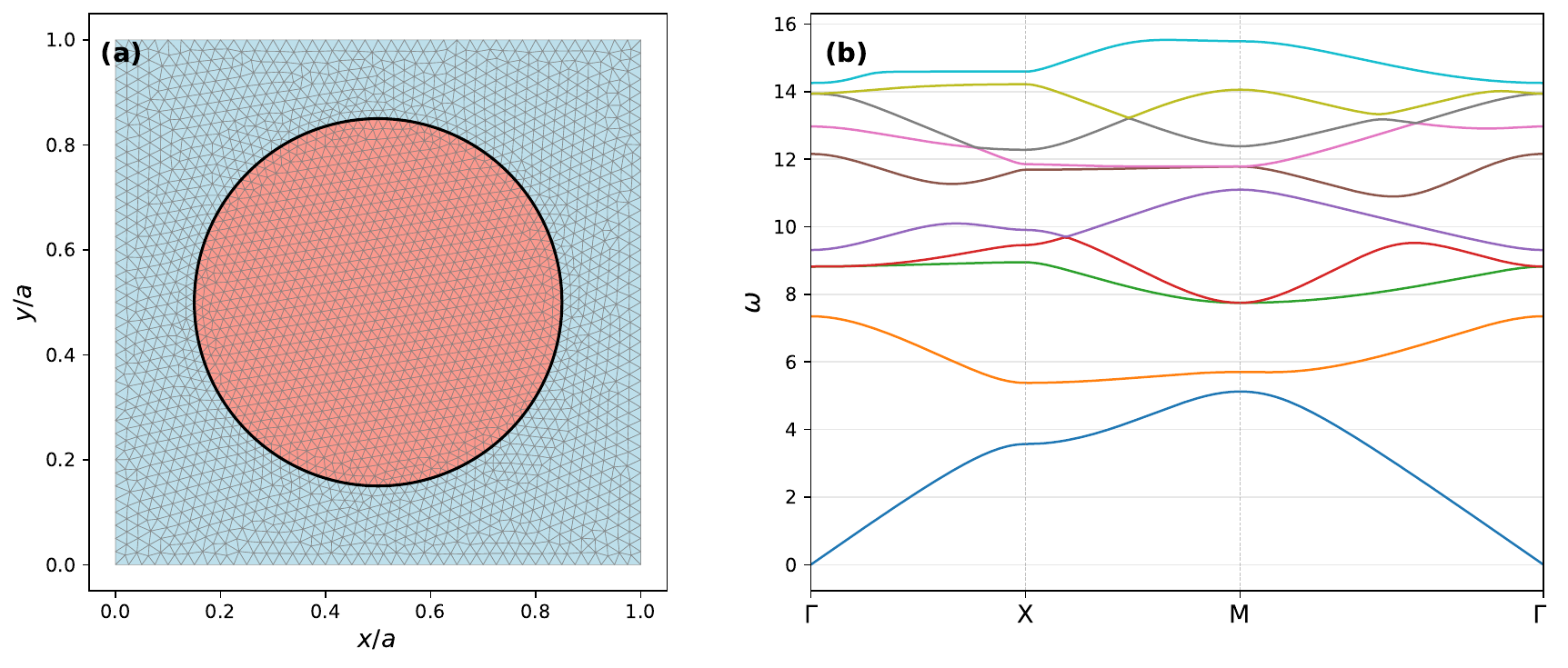}
  \caption{Two-dimensional phononic crystal with a circular inclusion
  ($r/a = 0.35$, $E_{\mathrm{inc}}/E_{\mathrm{mat}} = 12$).
  \textbf{(a)}~Triangular mesh of the unit cell; the matrix region is shown
  in blue and the inclusion in red.
  \textbf{(b)}~Band structure along the irreducible Brillouin zone boundary
  $\Gamma$--$X$--$M$--$\Gamma$, showing the first ten dispersion branches
  and two complete band gaps.}
  \label{fig:2d-mesh-bands}
\end{figure}

The band structure is computed along the boundary of the irreducible
Brillouin zone for the square lattice:
$\Gamma(0,0) \to X(\pi/a,0) \to M(\pi/a,\pi/a) \to \Gamma(0,0)$.
\Cref{fig:2d-mesh-bands}(b) shows the first ten dispersion branches. The
acoustic branch rises from $\omega = 0$ at $\Gamma$ and reaches its maximum
near the $M$ point. The higher bands exhibit
multiple near-crossings and narrow avoided crossings throughout the zone
boundary.

\subsection{SVD and the role of parameter-space dimension}
\label{sec:2d:svd}

The central prediction of the $n$-width theory developed in
\Cref{sec:analysis} is that the decay rate depends on the dimension $d$ of
the parameter space through the exponent $n^{1/d}$:
\begin{equation}
  \label{eq:2d-nwidth-recall}
  d_n(\calM_j, V) \leq C\,e^{-\beta\,n^{1/d}}.
\end{equation}
For the one-dimensional crystal of \Cref{sec:1d}, the Brillouin zone is an
interval ($d = 1$) and the bound reduces to a pure exponential
$e^{-\beta n}$, which was confirmed by the singular value decay in
\Cref{fig:phononic_1d}(c). For the present two-dimensional crystal, the
full Brillouin zone is a two-dimensional domain ($d = 2$), and the bound
becomes the stretched exponential $e^{-\beta\sqrt{n}}$. However, this
prediction applies only when the snapshots are sampled over a
two-dimensional region of $\kk$-space. If the snapshots are instead sampled
along a one-dimensional path---such as the IBZ boundary
$\Gamma$--$X$--$M$--$\Gamma$---the effective parameter dimension is $d = 1$
and pure exponential decay should be recovered, regardless of the ambient
spatial dimension of the crystal.

To test this, we perform two separate SVD analyses using the first $J = 3$
bands. In the first, the snapshot matrix is assembled from eigenvectors
sampled along the one-dimensional IBZ boundary path. In the second, the
snapshots are sampled over the interior of the two-dimensional irreducible
Brillouin zone---the triangle with vertices $\Gamma$, $X$, and $M$---using
a quasi-uniform grid. For each dataset, the singular values $\sigma_n$ of
the snapshot matrix are computed and normalized by $\sigma_1$.

\begin{figure}[htbp]
  \centering
  \includegraphics[width=\textwidth]{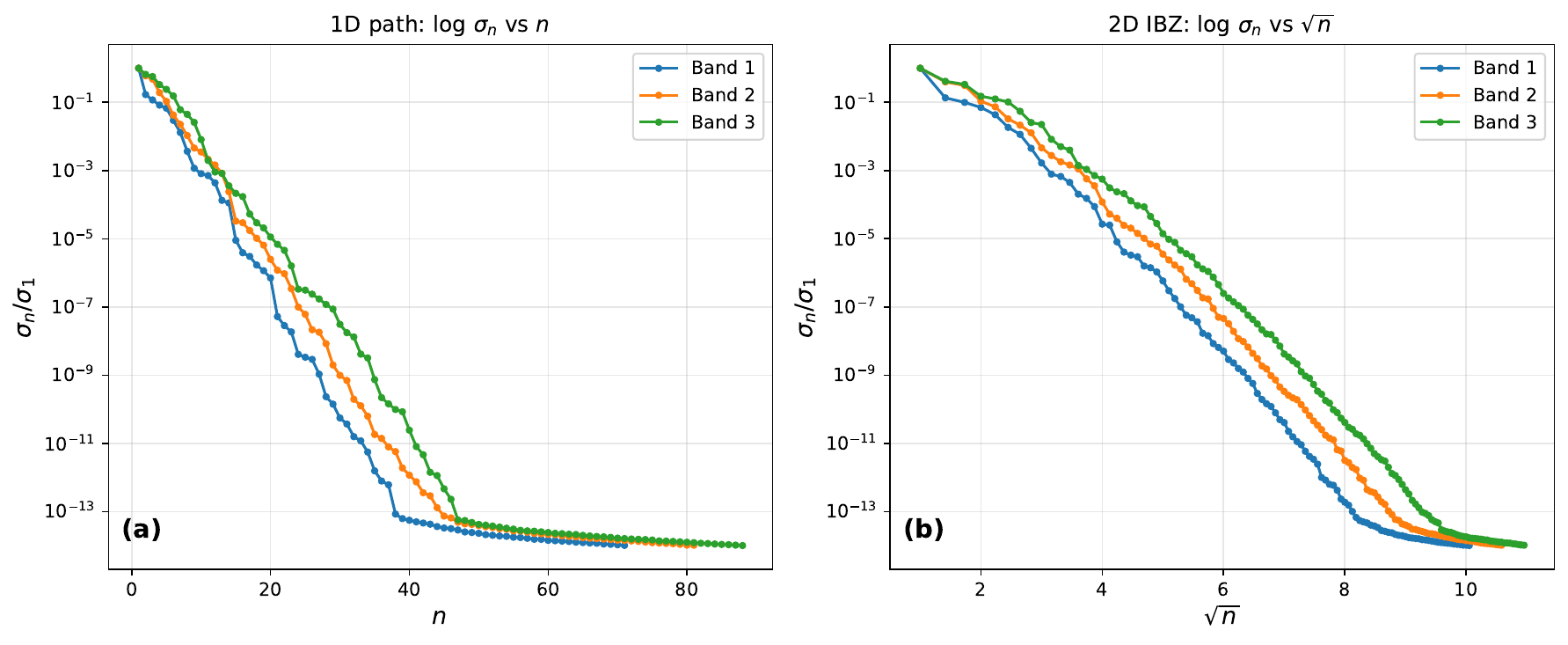}
  \caption{SVD decay of the snapshot matrices for the first three bands of
  the 2D phononic crystal.
  \textbf{(a)}~Snapshots sampled along the 1D IBZ boundary path: normalized
  singular values $\sigma_n/\sigma_1$ versus~$n$ on a semilog scale. The
  approximately linear decay confirms a pure exponential rate
  $e^{-\beta n}$, consistent with $d = 1$.
  \textbf{(b)}~Snapshots sampled over the 2D IBZ interior: $\sigma_n/\sigma_1$
  versus~$\sqrt{n}$. The approximately linear decay confirms the stretched
  exponential rate $e^{-\beta\sqrt{n}}$ predicted for $d = 2$.
  \Cref{thm:isolated-band}.}
  \label{fig:2d-svd}
\end{figure}

\Cref{fig:2d-svd}(a) shows $\log(\sigma_n/\sigma_1)$ versus~$n$ for the
three bands sampled along the one-dimensional IBZ boundary. All three curves
are approximately linear over the full range of significant singular values,
confirming pure exponential decay $\sigma_n \sim e^{-\beta n}$ consistent
with $d = 1$. Band~1, which has the largest spectral gap, decays fastest;
bands~2 and~3, with progressively narrower gaps, decay more slowly, in
agreement with the gap-controlled rate $\beta$ predicted by
\Cref{thm:isolated-band}.

\Cref{fig:2d-svd}(b) shows $\log(\sigma_n/\sigma_1)$ versus~$\sqrt{n}$
for the same three bands, now sampled over the two-dimensional IBZ. The
curves are approximately linear in $\sqrt{n}$, consistent with the
stretched exponential decay $\sigma_n \sim e^{-\beta\sqrt{n}}$ predicted by
the bound \eqref{eq:nwidth-isolated} with $d = 2$. The same data plotted
against~$n$ (not shown) produces concave curves, ruling out pure
exponential decay and confirming that the slower rate is genuine rather than
an artifact of insufficient sampling.

The practical consequence of the dimension-dependent exponent is visible in
the data. Along the one-dimensional IBZ boundary path, all three bands
reach a normalized singular value of $10^{-10}$ by approximately
$n \approx 30$--$40$ basis vectors. Over the two-dimensional IBZ interior,
achieving the same tolerance requires roughly $n \approx 50$--$70$ basis
vectors---a factor of roughly two increase in basis size. The increase from $d = 1$ to $d = 2$ is modest in absolute
terms---a few dozen additional vectors---but the distinction becomes more
significant at higher accuracy or for three-dimensional Brillouin zones,
where $d = 3$. For practical reduced-order
modeling, this means that computations targeting the full Brillouin zone of
a two- or three-dimensional crystal will require more basis vectors than
those targeting a one-dimensional path through it.

\section{Conclusions}
\label{sec:discussion}

This paper has established a rigorous optimality framework for reduced-order
methods used in phononic, acoustic, and photonic band structure calculations,
built on the theory of Kolmogorov $n$-widths.

The central result is that the Bloch eigenvalue problem is a particularly
clean instance of parametric holomorphy. The Bloch-transformed operators
$K(\kk)$ and $M(\kk)$ are entire holomorphic functions of the wave
vector~$\kk$, a consequence of the affine decomposition in the phase factors
$e^{i\kk\cdot\mathbf{a}_j}$. By Kato's analytic perturbation theory, the
eigenpairs inherit this holomorphy wherever the spectral gap is positive, with
a holomorphy radius bounded below by $\rho_j(\kk_0) \geq \delta_j(\kk_0) /
(2L)$. Applying classical results from approximation theory, the Kolmogorov
$n$-width of the solution manifold for an isolated band decays exponentially:
$d_n(\calM_j, V) \leq C\,e^{-\beta\,n^{1/d}}$, where the rate $\beta$ is
controlled by the minimum spectral gap between the band of interest and its
neighbors.

For multi-band computations, we showed that working with spectral projectors
rather than individual eigenvectors renders all internal
crossings---avoided, symmetry-enforced, or conical---irrelevant to the
convergence rate. Only the gap separating the cluster of $J$ targeted bands
from the remaining spectrum matters, and the $n$-width satisfies
$d_n(\calM_J, V) \leq C\,e^{-\beta\,(n-J)^{1/d}}$ for $n \geq J$. This
provides a unified treatment of all forms of spectral degeneracy encountered
in practice.

Numerical experiments on a one-dimensional phononic crystal confirmed the
theoretical predictions. The singular values of the snapshot matrix---an
upper bound on the $n$-width---decayed exponentially for all bands, with
rates consistent with the spectral gap analysis. Both an oracle greedy
algorithm and a practical residual-based greedy tracked the SVD-optimal
subspace within a small constant factor, reaching machine precision with
approximately $J + 20$ basis vectors for the first ten bands. The greedy
selection sequence revealed that the most informative sampling locations are
the Brillouin zone boundaries, where eigenfunctions differ most from those
at the zone center, and that the highest bands---those with the smallest
spectral gaps---are selected first.

The two-dimensional numerical experiments of \Cref{sec:2d} confirmed the
role of the parameter-space dimension in the $n$-width bound. By comparing
singular value decay for snapshots sampled along the one-dimensional IBZ
boundary path against snapshots sampled over the two-dimensional IBZ
interior, we verified that the exponent $n^{1/d}$ in the bound
$d_n \leq C\,e^{-\beta\,n^{1/d}}$ is genuine: the one-dimensional path
produces pure exponential decay in~$n$, while the two-dimensional domain
produces the stretched exponential decay in~$\sqrt{n}$ predicted by the
theory. The basis size required for a given tolerance approximately doubled
from $d = 1$ to $d = 2$, consistent with the scaling
$n(\varepsilon) \sim |\log\varepsilon|^d$. 

These results have several implications for the design and evaluation of
existing methods. First, they provide a sharp lower bound against which any
reduced-order method can be measured: a method whose error decays
exponentially with basis size at a rate comparable to $\beta$ is
near-optimal and cannot be substantially improved by any other linear
approach. Second, the analysis offers a principled justification for the
empirical basis selection strategies used in methods such as
RBME~\cite{Hussein2009}: sampling eigenvectors at high-symmetry points of
the Brillouin zone is effective precisely because these points are near the
zone boundaries where the eigenvectors are most poorly represented by an
interior-point initialization---a fact independently discovered by the
greedy algorithm. Third, the spectral projector formulation explains why
methods that retain all eigenvectors within a frequency window, such as
BMS~\cite{KrattigerHussein2014}, exhibit robust convergence even in the
presence of band crossings: they implicitly approximate the spectral
subspace rather than individual eigenvectors, and the relevant gap is the
cluster gap at the window boundary.

Several directions for future work suggest themselves. The present numerical
experiments are limited to one spatial dimension, where the Brillouin zone
is an interval and all eigenvalues are generically simple. Extending the
computations to two- and three-dimensional crystals would test the predicted
$e^{-\beta\,n^{1/d}}$ decay, where the stretched exponent $n^{1/d}$
reflects the higher-dimensional parameter space, and would probe the
practical impact of symmetry-enforced degeneracies and conical
intersections. A second direction concerns the tightness of the bound: while
the exponential rate $\beta \sim \delta_j^*/(2L)$ correctly identifies the
spectral gap as the controlling quantity, it may be conservative, and
sharper estimates exploiting the specific structure of the Bloch
problem---such as the periodicity of the Brillouin zone and the finite
number of terms in the affine decomposition---could yield improved rates. A
third direction is the extension of this analysis to nonlinear parameter
dependence, as arises in the inverse band structure problem
$\kk(\omega)$~\cite{PalermoMarzani2016}, where the roles of frequency and
wave vector are exchanged and the parametric structure differs. Finally, the
framework developed here applies equally to electronic band structure
calculations, where the connection to Wannier function
theory~\cite{MarzariVanderbilt1997,MarzariMostofiYatesSouzaVanderbilt2012}
and the role of topological obstructions in limiting smooth frame
construction merit further investigation within the $n$-width setting.

\appendix
\section{Example: 1D Homogeneous Bar with Two Linear Elements}
\label{app:1d-example}

Consider a one-dimensional phononic crystal consisting of a homogeneous bar
with Young's modulus $E$, mass density $\rho$, and cross-sectional area $A$.
The unit cell is the interval $\Omega = [0, a]$, and the lattice vector is
simply $\mathbf{a}_1 = a$. We discretize the unit cell with two linear
(two-node) bar elements of equal length $h = a/2$, placing nodes at
$x_0 = 0$, $x_1 = a/2$, and $x_2 = a$. Assembling the two elements over the three-node mesh yields the $3 \times 3$
global matrices
\begin{equation}
  \label{eq:app-global}
  K_{\mathrm{g}} = \frac{EA}{h}
  \begin{pmatrix} 1 & -1 & 0 \\ -1 & 2 & -1 \\ 0 & -1 & 1 \end{pmatrix},
  \qquad
  M_{\mathrm{g}} = \frac{\rho A h}{6}
  \begin{pmatrix} 2 & 1 & 0 \\ 1 & 4 & 1 \\ 0 & 1 & 2 \end{pmatrix},
\end{equation}
with rows and columns ordered as $(u_0, u_1, u_2)$. The Bloch-periodic boundary condition requires
$u_2 = e^{ika}\,u_0$, relating the displacement at the right boundary of the
unit cell to that at the left boundary through the phase factor
$z = e^{ika}$. This constraint eliminates one degree of freedom, reducing the
system from three unknowns to $N = 2$ independent degrees of freedom
$\mathbf{q} = (u_0,\; u_1)^\top$. The constraint is enforced through the transformation
$\mathbf{u}_{\mathrm{g}} = T(k)\,\mathbf{q}$, where
\begin{equation}
  \label{eq:app-transform}
  T(k) = \begin{pmatrix} 1 & 0 \\ 0 & 1 \\ e^{ika} & 0 \end{pmatrix}.
\end{equation}
The Bloch-reduced stiffness and mass matrices are obtained by
congruence:
\begin{equation}
  \label{eq:app-bloch-assembly}
  \hat{K}(k) = T(k)^H\,K_{\mathrm{g}}\,T(k), \qquad
  \hat{M}(k) = T(k)^H\,M_{\mathrm{g}}\,T(k),
\end{equation}
where the superscript $H$ denotes the conjugate transpose, required because
$T(k)$ is complex-valued. Carrying out the matrix multiplications yields the $2 \times 2$
Bloch operators
\begin{equation}
  \label{eq:app-bloch-K}
  \hat{K}(k) = \frac{2EA}{a}
  \begin{pmatrix}
    2 & -1 - e^{-ika} \\
    -1 - e^{ika} & 2
  \end{pmatrix},
\end{equation}
\begin{equation}
  \label{eq:app-bloch-M}
  \hat{M}(k) = \frac{\rho A a}{12}
  \begin{pmatrix}
    4 & 1 + e^{-ika} \\
    1 + e^{ika} & 4
  \end{pmatrix}.
\end{equation}

\subsection{Affine decomposition}
\label{app:1d:affine}

The $k$-dependence of the Bloch matrices \eqref{eq:app-bloch-K}--\eqref{eq:app-bloch-M} is
entirely contained in the phase factors $e^{\pm ika}$, providing an explicit
instance of the affine decomposition \eqref{eq:affine-decomp}. For the
stiffness matrix:
\begin{equation}
  \label{eq:app-affine-K}
  \hat{K}(k) = K_0 + e^{ika}\,K_1 + e^{-ika}\,K_2,
\end{equation}
where the $k$-independent component matrices are
\begin{equation}
  \label{eq:app-K-components}
  K_0 = \frac{2EA}{a}
  \begin{pmatrix} 2 & -1 \\ -1 & 2 \end{pmatrix}, \quad
  K_1 = \frac{2EA}{a}
  \begin{pmatrix} 0 & 0 \\ -1 & 0 \end{pmatrix}, \quad
  K_2 = \frac{2EA}{a}
  \begin{pmatrix} 0 & -1 \\ 0 & 0 \end{pmatrix}.
\end{equation}
The mass matrix admits the analogous decomposition
$\hat{M}(k) = M_0 + e^{ika}\,M_1 + e^{-ika}\,M_2$ with
\begin{equation}
  \label{eq:app-M-components}
  M_0 = \frac{\rho Aa}{12}
  \begin{pmatrix} 4 & 1 \\ 1 & 4 \end{pmatrix}, \quad
  M_1 = \frac{\rho Aa}{12}
  \begin{pmatrix} 0 & 0 \\ 1 & 0 \end{pmatrix}, \quad
  M_2 = \frac{\rho Aa}{12}
  \begin{pmatrix} 0 & 1 \\ 0 & 0 \end{pmatrix}.
\end{equation}
Note that $K_2 = K_1^H$ and $M_2 = M_1^H$, reflecting the Hermitian symmetry
$\hat{K}(k)^H = \hat{K}(k)$ and $\hat{M}(k)^H = \hat{M}(k)$, which holds
for all $k \in \R$. The number of terms in the decomposition ($Q = 2$,
corresponding to the phase factors $e^{\pm ika}$) is determined by the mesh
connectivity across the unit cell boundary---here, only nearest-neighbor
coupling across the boundary exists---and does not depend on the mesh
refinement within the cell.

\subsection{Eigenvalues and eigenvectors}
\label{app:1d:eigenpairs}

The generalized eigenvalue problem
$\hat{K}(k)\,\mathbf{q} = \omega^2\,\hat{M}(k)\,\mathbf{q}$ for the
$2 \times 2$ system \eqref{eq:app-bloch-K}--\eqref{eq:app-bloch-M} can be
solved in closed form. The eigenvectors of this system are
\begin{equation}
  \label{eq:app-evecs}
  \mathbf{u}_1(k) = \begin{pmatrix} 1 \\ e^{ika/2} \end{pmatrix}, \qquad
  \mathbf{u}_2(k) = \begin{pmatrix} 1 \\ -e^{ika/2} \end{pmatrix},
\end{equation}
which can be verified by direct substitution. The corresponding eigenvalues
are
\begin{equation}
  \label{eq:app-evals}
  \omega_1^2(k) = \frac{24\,c_0^2}{a^2}\cdot
  \frac{1 - \cos(ka/2)}{2 + \cos(ka/2)}, \qquad
  \omega_2^2(k) = \frac{24\,c_0^2}{a^2}\cdot
  \frac{1 + \cos(ka/2)}{2 - \cos(ka/2)},
\end{equation}
where $c_0 = \sqrt{E/\rho}$ is the bar wave speed.

\subsection{The solution manifold}
\label{app:1d:manifold}

The solution manifold for the first band is
\begin{equation}
  \label{eq:app-manifold}
  \calM_1 = \bigl\{\mathbf{u}_1(k) \in \C^2 :
  k \in [-\pi/a,\;\pi/a]\bigr\}
  = \Bigl\{\begin{pmatrix} 1 \\ e^{ika/2} \end{pmatrix} :
  k \in [-\pi/a,\;\pi/a]\Bigr\}.
\end{equation}
The first component is constant, while the second traces a unit-circle arc in
$\C$ as $k$ sweeps the Brillouin zone. Separating real and imaginary parts
and viewing $\C^2 \cong \R^4$, the manifold becomes
\begin{equation}
  \label{eq:app-manifold-real}
  \calM_1 = \bigl\{
  \bigl(1,\; 0,\; \cos(ka/2),\; \sin(ka/2)\bigr)^\top :
  k \in [-\pi/a,\;\pi/a]
  \bigr\},
\end{equation}
which is a circular arc lying in
$\R^4$. To represent this manifold exactly, one needs a three-dimensional linear subspace: two dimensions span the circular arc traced by the third and fourth components, and a third is required because the constant first component shifts the manifold away from the origin, which every linear subspace must contain.

\bibliographystyle{unsrt}
\bibliography{references}

\end{document}